\documentclass[english,twoside,11pt]{article}

\usepackage[usenames]{color}
\usepackage[german,english]{babel}
\usepackage{amsmath}
\usepackage{amssymb,latexsym}
\usepackage{moreverb}
\usepackage{url}
\usepackage{graphicx}
\usepackage{longtable}
\usepackage{rotating}
\usepackage{lscape}
\usepackage{booktabs}
\usepackage{tabularx}
\usepackage{psfrag}

\usepackage{multicol}
\usepackage{gastex}
\newtheorem{deff}{Definition :}
\newtheorem{lem}[deff]{Lemma :}

\newtheorem{thm}[deff]{Theorem :}

\newtheorem{cor}[deff]{Corollary :}

\newtheorem{example}[deff]{Example :}

\newtheorem{claim}[deff]{Claim :}

\makeatletter

\title{Some Semi - Equivelar Maps}
\author{Ashish K. Upadhyay,  Anand K. Tiwari and Dipendu Maity\\[1mm]
Department of Mathematics\\Indian Institute of Technology Patna\\
Patliputra Colony, Patna -- 800\,013, India\\
{\small \{upadhyay, anand, dipendumaity\}@iitp.ac.in}}

\date{\today}

\begin{document}

\selectlanguage{english}

\maketitle

\hrule

\begin{abstract}

Semi-Equivelar maps are generalizations of Archimedean Solids (as are equivelar maps of the Platonic solids) to the surfaces other than $2-$Sphere. We classify some semi equivelar maps on surface of Euler characteristic $-1$ and show that none of these are vertex transitive. We establish existence of $12$-covered triangulations for this surface. We further construct double cover of these maps to show existence of semi-equivelar maps on the surface of double torus. We also construct several semi-equivelar maps on the surfaces of Euler characteristics $-8$ and $-10$ and on non-orientable surface of Euler characteristics $-2$.


\end{abstract}

\hrule

\bigskip

{\bf AMS Subject Classification\,:} 57M10, 57M20, 52B50, 52B70, 52C20

{\bf Keywords\,:} Semi-Equivelar Maps, $d$-Covered Triangulations, Equivelar Triangulations.


\section{Introduction and results}

As is well known, equivelar triangulations, also known as degree regular triangulations of surfaces and in more generality equivelar maps on surfaces are generalizations of the maps on surfaces of five Platonic solids to the surfaces other than sphere. Here we attempt to study generalizations of Archimedean solids to some surfaces of negative Euler characteristics. Their study for the surfaces of non-negative Euler characteristics has been carried out in \cite{upadhyay&tiwari2}. We call such objects as {\em Semi-Equivelar Maps} or SEM(s).

The following definitions are given in \cite{datta} and we reproduce it here for the sake of completeness and ready reference. A {\em p-Cycle}, denoted $C_p$,  is a finite connected 2-regular graph with $p$ vertices. A 2-dimensional {\em Polyhedral Complex} $K$ is a collection of $p_{i}$-cycles, where $\{p_{i}\colon 1 \leq i \leq n\}$ is a set of positive integers $\geq 3$, together with vertices and edges in the cycles such that the intersection of any two cycles is empty, a vertex or is an edge. The cycles are called faces of $K$. The symbols $V(K)$ and $EG(K)$ respectively denote the set of vertices and edges of $K$. A polyhedral complex $K$ is called a {\em Polyhedral 2-manifold} if for each vertex $v$ the faces containing $v$ are of the form $C_{p_1}, \ldots, C_{p_m}$ where $C_{p_1}\cap C_{p_2}, \ldots C_{p_{m-1}}\cap C_{p_m}$,  and  $C_{p_m}\cap C_{p_1}$ are edges for some $m \geq 3$. A connected finite polyhedral 2-manifold is called a {\em Polyhedral Map}. We will use the term {\em map} for a polyhedral map.

We associate a geometric object $|K|$ to a polyhedral complex $K$ as follows\,: corresponding to each $p$-cycle $C_p$ in $K$, consider a $p$-gon $D_p$ whose boundary cycle is $C_p$. Then $|K|$ is union of all such $p$-gons and is called the {\em geometric carrier} of $K$. The complex $K$ is said to be connected (resp. orientable) if $|K|$ is connected (resp. orientable) topological space. Between any two polyhedral complexes $K_1$ and $K_2$ we define an isomorphism to be a map $f\colon K_1 \longrightarrow K_2$ such that $f|_{V(K_1)}\colon V(K_1) \longrightarrow V(K_2)$ is a bijection and $f(\sigma)$ is a cell in $K_2$ if and only if $\sigma$ is a cell in $K_1$. If $K_1 = K_2$ then $f$ is called an automorphism of $K_1$. The set of all automorphisms of a polyhedral complex $K$ form a group under the operation composition of maps. This group is called the group of automorphisms of $K$. If this group acts transitively on the set $V(K)$ then the complex is called a {\em vertex transitive} complex. Some vertex transitive maps of Euler characteristic 0 have been studied in \cite{Babai}.


The {\em Face sequence} of a vertex $v$ in a map, see figure in Example \ref{e2}, is a finite sequence $(a^p, b^q, ...., m^r)$ of powers of positive integers  $a, b, ..., m \geq 3$ and $p, q, ..., r \geq 1$, such that through the vertex $v$, $p$ numbers of $C_a$, $q$ numbers of $C_b$, $\ldots$, $r$ numbers of $C_m$ are incidents. A map $K$ is said to be {\em Semi-Equivelar} if face sequence of each vertex of $K$ is same. Thus, for example the face sequence  of a vertex in the maps of Example $\ref{e1}$ is $(3^5, 4)$. In \cite{upadhyay&tiwari2}, maps with face sequence $(3^3, 4^2)$ and $(3^2, 4, 3, 4)$ have been considered.

A triangulation of a connected closed surface is called {\em Equivelar}  (or degree regular) if each of its vertices have  the same degree. Thus a $d-$equivelar triangulation is a SEM of type $(3^d)$. A triangulation is called {\em $d$-covered} if each edge of the triangulation is incident with a vertex of degree $d$. In articles \cite{dattaupadhyay1} and \cite{dattaupadhyay2} equivelar triangulations have been studied for Euler characteristics  0 and -2. In \cite{NegamiNak} Negami and Nakamoto studied $d$-covered triangulations and asked the question about existence of such triangulations on a surface of Euler charateristic $\chi$ with the condition that $d = 2\,\lfloor\displaystyle\frac{5 + \sqrt{49 - 24\chi}}{2}\rfloor$ see also \cite{aku1}. In the articles \cite{LSTU} questions about existence of such triangulations for $-2\leq \chi \leq -127$ and whenever $n = \displaystyle\frac{\chi\,d}{d - 6}$ is an integer is considered. There, use of equivelar triangulations of surfaces to construct the required $d$-covered triangulations has been made. It is well known that the equivelar triangulations do not exist for surface of Euler characteristic $-1$ and some study about maps on this surface has been made in \cite{AltBrehm}. The current work is motivated by an attempt to search for existence of $12$-covered triangulations on the surface of Euler characteristic $-1$. We answer the question in affirmative. We construct and classify some semi equivelar maps on this surface. A glance at the Euler's formula ${\chi = \,{no. \,of\, vertices}}$ - ${no. \,of\, edges}$ + ${no.\, of\, faces}$ , indicates that at each vertex, if we allow one of the faces to be a quadrangle and choose smallest possible number of triangles in such a way that curvature is negative then it might be possible to have some maps on this surface. Such a map is an example of what we defined to be a {\em Semi-Equivelar map} on a surface. In current work, we allowed one square and five triangles at each vertex to discover  that if $N$ denotes number of vertices, the Euler's equation gives $\chi = N (\frac{-1}{12})$. Thus if we take $N = 12$ then  we may obtain a SEM of type $(3^5, 4)$ on the non - orientable surface of Euler characteristic $-1$. In the following lines we will be defining a procedure to add handles to a SEM. This process is not new and appeared earlier in $\cite{McSchWills}$.

Let $C_l(v_1, v_2, \ldots, v_l)$ and $Z_l(u_1, u_2, \ldots, u_l)$ denote  cycles of length $l$. We define a cylinder $C_{ll}(C_l, Z_l)$ with boundary components $C_l$ and $Z_l$ to be the complex with vertex set $\{v_1, v_2, \ldots, v_l, u_1$, $u_2, \ldots, u_l\}$ and facets $\{v_1v_2u_{i_2}u_{i_1}, v_2v_3u_{i_3}u_{i_2}, \ldots , v_1v_lu_{i_l}u_{i_1}\}$. If $K$ denotes a SEM of type $(3^5, 4)$, it is possible to add a cylinder to $K$ to obtain another map as follows: let $Q_1$ and $Q_2$ denote two quadrangular facets of $K$ such that $V(Q_1)\cap V(Q_2) = \empty $. Remove the interior of $Q_1$ and $Q_2$ to obtain two disjoint cycles $\partial(Q_1)$ and $\partial(Q_2)$ of length four as boundary components of $K$. One may now plug in a cylinder $C_{44}(\partial(Q_1), \partial(Q_2) )$ by identifying the boundaries. If one is careful in choosing pairs $(Q_1, Q_2)$ so as to preserve a semi-equivelar type for $K$, we do obtain a genuine map from $K$. One can see that by this the Euler characteristic is increased by $-2$. It is sometimes needed to triangulate the cylinder suitably (so as to retain the semi-equivelar type) and obtain an addition of triangular facets instead of quadrangular faces. We perform two types of cylinder addition $(i)$ we adjoin a cylinder to the boundary components of $K\setminus\{Q_1, Q_2\}$, or $(ii)$ we adjoin a cylinder in the boundary of $K\setminus Q_1 \bigcup L\setminus Q_2$ for two SEMs $K, L$ such that $Q_1 \in K$ and $Q_2 \in L$.  We use both these types of cylinder additions to obtain SEM of types $\{3^5, 4^2\}$ and $\{3^7, 4\}$ on surfaces of Euler characteristics $-8$ and $-10$, see Examples \ref{e4} and \ref{e5}.

Let $EG(K)$ be the edge graph of a map $K$ and $V(K) = \{v_1, v_2, \ldots, v_n\}$. Let $L_{K}(v_i) = \{v_j \in V(K) \colon v_iv_j \in EG(K)\}$. For $0 \leq t \leq n$ we define a graph $G_{t}(K)$ with $V(G_t(K)) = V(K)$ and $v_iv_j \in EG(G_t(K))$ if $|L_{K}(v_i) \bigcap L_{K}(v_j)| = t$, in other words the number of elements in  the set $L_{K}(v_i) \bigcap L_{K}(v_j)$ is $t$. This graph was introduced in \cite{dattaupadhyay2} by B. Datta. Moreover if $K$ and $K'$ are two isomorphic maps then $G_i(K) \cong G_i(K')$ for each $i$. For computations, we have used {\tt GAP} \cite{gap}. We have also computed reduced homology groups of the objects using \cite{hecken}. We classify all the semi-equivelar maps of type $(3^5, 4)$ on the surface of Euler characteristics -1 and show that there are precisely three such objects up to isomorphism. In other words we show that \,:


\begin{thm}\label{t1} If $K$ is a semi equivelar maps of type $(3^5, 4)$ on the surface of Euler characteristic $-1$ then $K$ is isomorphic to one of $K_1$, $K_2$ or $K_3$ given in Example \ref{e1}.
\end{thm}


\begin{cor} There exists a $12$-covered triangulation of the surface of Euler characteristics $-1$.
\end{cor}
\begin{proof} To each face of the map add the barycenter and join each vertices of the face with this newly introduced vertex. This process is called stacking a face, see \cite{LSTU}. The proof now follows by stacking each face of the $(3^5, 4)$ SEM on the surfaces of Euler characteristic $-1$. The resulting triangulation is $12$-covered.
\end{proof}$\hfill\Box$


In $\cite{karabasNedela}$, Karabas and Nedela have presented a census of vertex transitive Archimedean solids of genus two. This census includes one SEM of type $(3^5, 4)$ on the surface of double torus on 24 vertices. Here we construct more such maps on the orientable surface of genus two as double covers of the SEMs of same type on surface of Euler characteristic $-1$. For each of these maps we present their double covers in Example 2 which turn out to be mutually non - isomorphic. We prove that\,:


\begin{thm}\label{t2} There exist at least four SEM of type $(3^5, 4)$ on the surface of Euler characteristic $-2$. Three of these is orientable and one is non orientable. None of these maps are vertex transitive.
\end{thm}

\begin{cor} There are at least five SEMs of type $(3^5, 4)$ on the surface of Euler characteristic $-2$. Four of these are orientable and one is non-orientable. Among the orientable SEMs one is vertex transitive and the remaining are not. $\hfill\Box$
\end{cor}

We further show the following\,:

\begin{thm}\label{t3} There exist at least 10 SEMs of types $(3^5, 4^2)$ on the surface of Euler characteristic $-8$. Two of these are orientable and remaining are non - orientable.
\end{thm}


\begin{thm}\label{t4} There exist at least 11 SEMs of types $(3^7, 4)$ on the surface of Euler characteristic $-10$. Two of these are orientable and remaining are non - orientable.
\end{thm}



\section{Examples}

\begin{example}\label{e1} {\bf Some Semi Equivelar Maps on surface of Euler Characteristics -1\,:}

\bigskip

\noindent{$\bf K_1$} = \{012, 017, 045, 056, 067, 128, 158, 15u, 236, 267, 278, 34v, 369, 39u, 3uv, 45u, 49u, 49v, 78v, 89v, 0234, 17vu, 5698\}

\bigskip

\noindent{$\bf K_2$} = \{012, 017, 045, 056, 067, 129, 17v, 189, 238, 268, 269, 34v, 389, 39u, 3uv, 45u, 47u, 47v, 568, 5uv, 0234, 185v, 67u9\}

\bigskip

\noindent{$\bf K_3$} = \{012, 017, 045, 056, 067, 129, 178, 19v, 238, 268, 269, 34v, 378, 37u, 3uv, 45u, 49u, 49v, 568, 5uv, 0234, 85v1, 67u9\}

The Graphs $EG(G_{6}(K_1)) = \emptyset$, $EG(G_2(K_1) ) =\{ [ 2, 4 ], [ 7, 10 ] \}$. $EG(G_2(K_2)) = \{ [ 2, 4 ], [ 3, 12 ] \}$ and $EG(G_6(K_2)) = \{[ 1, 6 ], [ 5, 7 ] \}$. Also, $EG(G_2(K_3)) = \emptyset$ and $EG(G_6(K_3)) = \{[ 1, 6 ], [ 8, 12 ] \}$. Therefore, $K_1 \not\cong K_2$, $K_1 \not\cong K_3$ and $K_2 \not\cong K_3$. A look at these graphs one can easily deduce that $K_1$, $K_2$ and $K_3$ are not vertex transitive.
\end{example}


\hrule


\begin{example}\label{e2} {\bf $(a)\,$ Some face sequences of vertex $v = 0$\:}

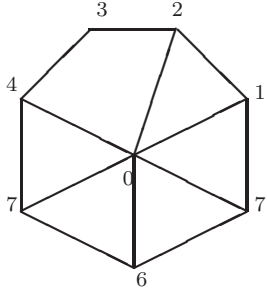
\setlength{\unitlength}{5mm}
\begin{picture}(0,0)(0,4)

\thicklines
\put(4,-2.5){\line(2,1){3}}
\put(7,-1){\line(0,1){3}}
\put(1,2){\line(1,1){1.82}}
\put(5.16,3.83){\line(1,-1){1.85}}
\put(1,-1){\line(0,1){3}}
\put(1,-1){\line(2,-1){3}}
\put(2.8,3.85){\line(1,0){2.29}}

\put(4,0.5){\line(2,-1){3}}
\put(4,0.5){\line(2,1){3}}
\put(1,2){\line(2,-1){3}}
\put(4.,-2.5){\line(0,1){3}}
\put(1,-1){\line(2,1){3}}
\put(4.,0.5){\line(1,3){1.12}}

\put(3.7,-0.3){\scriptsize $0$}
\put(4.05,-3){\scriptsize $6$}
\put(7.2,-1){\scriptsize $7$}
\put(7.2,2){\scriptsize $1$}
\put(5,4.2){\scriptsize $2$}
\put(3,4.2){\scriptsize $3$}
\put(0.6,2.2){\scriptsize $4$}
\put(0.6,-1){\scriptsize $7$}

\put(3,-4){\scriptsize{$Type\,: (3^5, 4)$}}

\end{picture}

\setlength{\unitlength}{5mm}
\begin{picture}(0,0)(-10,3.5)

\thicklines

\put(4,-2.5){\line(2,1){3}}
\put(7,-1){\line(0,1){3}}
\put(1,2){\line(1,1){1.82}}
\put(5.16,3.83){\line(1,-1){1.85}}
\put(1,-1){\line(0,1){3}}
\put(1,-1){\line(2,-1){3}}
\put(2.85,3.84){\line(1,0){2.23}}

\put(4,0.5){\line(2,1){3}}
\put(1,2){\line(2,-1){3}}
\put(4.,-2.5){\line(0,1){3}}
\put(4.,0.5){\line(1,3){1.12}}
\put(2.9,3.8){\line(1,-3){1.1}}

\put(3.7,-0.3){\scriptsize $0$}
 \put(4.05,-3){\scriptsize $6$}
\put(7.2,-1){\scriptsize $7$}
\put(7.2,2){\scriptsize $1$}
\put(5,4.2){\scriptsize $2$}
\put(3,4.2){\scriptsize $3$}
\put(0.6,2.2){\scriptsize $4$}
\put(0.6,-1){\scriptsize $7$}

\put(3,-4){\scriptsize{$Type\,: (3^3, 4^2)$}}

\end{picture}\
\setlength{\unitlength}{5mm}
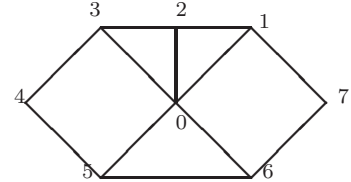
\begin{picture}(0,0)(-22,4)
\thicklines

\put(0,4){\line(1,0){4}}
\put(0,0){\line(1,0){4}}
\put(2,2){\line(0,1){2}}
\put(2,2){\line(1,1){2}}
\put(0,4){\line(1,-1){2}}
\put(-2,2){\line(1,1){2}}
\put(4,4){\line(1,-1){2}}
\put(-2,2){\line(1,-1){2}}
\put(4,0){\line(1,1){2}}
\put(0,0){\line(1,1){2}}
\put(2,2){\line(1,-1){2}}
\put(2,1.3){\scriptsize$0$}
\put(4.2,4){\scriptsize$1$}
\put(2,4.3){\scriptsize$2$}
\put(-.3,4.3){\scriptsize$3$}
\put(-2.3,2){\scriptsize$4$}
\put(-.5,0){\scriptsize$5$}
\put(4.3,0){\scriptsize$6$}
\put(6.3,2){\scriptsize$7$}
\put(2,-2.5){\scriptsize{$Type\,: (3^2, 4, 3, 4)$}}

\end{picture}

\vspace{1.4in}

\end{example}

\smallskip

\noindent {\bf $(b)\,$ Figure showing two types of Cylinders\,:}

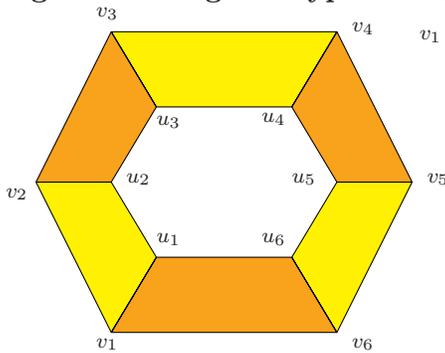
\begin{center}

\begin{picture}(120,40)(0,0)\nullfont
  \put(20,20){
    \unitlength=2mm
    \drawpolygon[fillcolor=Yellow,Nframe=n,arcradius=0](-10,-10)(-7,-5)(-10,0)(-15,0)
    \drawpolygon[fillcolor=YellowOrange,Nframe=n,arcradius=0](-15,0)(-10,0)(-7,5)(-10,10)
    \drawpolygon[fillcolor=Yellow,Nframe=n,arcradius=0](-10,10)(-7,5)(2,5)(5,10)
    \drawpolygon[fillcolor=YellowOrange,Nframe=n,arcradius=0](2,5)(5,10)(10,0)(5,0)
    \drawpolygon[fillcolor=Yellow,Nframe=n,arcradius=0](10,0)(5,0)(2,-5)(5,-10)
    \drawpolygon[fillcolor=YellowOrange,Nframe=n,arcradius=0](5,-10)(2,-5)(-7,-5)(-10,-10)
    \drawpolygon(-10,-10)(-15,0)(-10,10)(5,10,)(10,0)(5,-10)
     \drawpolygon(-7,-5)(-10,0)(-7,5)(2,5,)(5,0)(2,-5)
     \drawpolygon(-7,-5)(-10,-10)
     \drawpolygon(-10,0)(-15,0)
      \drawpolygon(5,0)(10,0)
      \drawpolygon(-7, 5)(-10, 10)
       \drawpolygon(2, 5)(5, 10)
        \drawpolygon(2,-5)(5,-10)
\put(-11,-11){\scriptsize $v_1$}
\put(-17,-1){\scriptsize $v_2$}
\put(-11,11){\scriptsize $v_3$}
\put(6,10){\scriptsize $v_4$}
\put(11,0){\scriptsize $v_5$}
\put(6,-11){\scriptsize $v_6$}
\put(-7,-4){\scriptsize $u_1$}
\put(-9,0){\scriptsize $u_2$}
\put(-7,3.8){\scriptsize $u_3$}
\put(0,4){\scriptsize $u_4$}
\put(2,0){\scriptsize $u_5$}
\put(0,-4){\scriptsize $u_6$}

\put(-5,-13){\scriptsize $C_{66}(v_1\dots v_6, u_1\dots u_6)$}
  }



  \put(85,19){
      \unitlength=2mm
    \drawpolygon[fillcolor=SpringGreen,Nframe=n,arcradius=0](-20,10)(-12,5)(20,10)
    \drawpolygon[fillcolor=YellowGreen,Nframe=n,arcradius=0](-12,5)(20,10)(5,5)
    \drawpolygon[fillcolor=SpringGreen,Nframe=n,arcradius=0](5,5)(20,10)(-10,-12)
    \drawpolygon[fillcolor=YellowGreen,Nframe=n,arcradius=0](5,5)(-10,-12)(-8,-5)
    \drawpolygon[fillcolor=SpringGreen,Nframe=n,arcradius=0](-10,-12)(-8,-5)(-20,10)
    \drawpolygon[fillcolor=YellowGreen,Nframe=n,arcradius=0](-20,10)(-12,5)(-8,-5)
    \drawpolygon(-20,10)(20,10)(-10,-12)
    \drawpolygon(-12,5)(5,5)(-8,-5)
     \drawline[AHnb=0](-12,5)(-20,10)
     \drawline[AHnb=0](20,10)(5,5)
      \drawline[AHnb=0](-10,-12)(-8,-5)
      \drawline[AHnb=0](20, 10)(-12, 5)
       \drawline[AHnb=0](-20,10)(-8, -5)
\drawline[AHnb=0](-10,-12)(5,5)
\put(-22,10){\scriptsize $v_1$}
\put(22,10){\scriptsize $v_2$}
\put(-8,-12){\scriptsize $v_3$}
\put(-11,3){\scriptsize $u_1$}
\put(0,3){\scriptsize $u_2$}
\put(-8,-3){\scriptsize $u_3$}

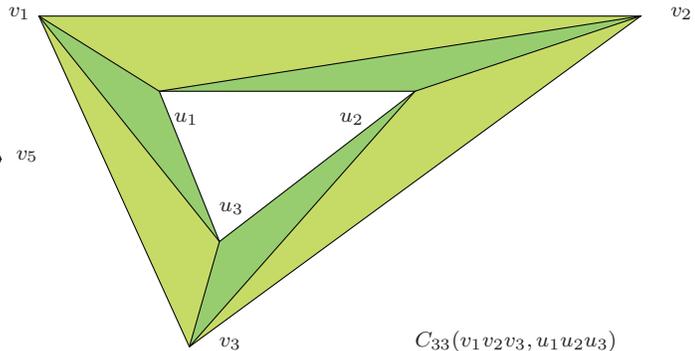
\put(5,-12){\scriptsize $C_{33}(v_1v_2v_3, u_1u_2u_3)$}
  }
\end{picture}

\vspace{0.2in}


\end{center}
\hrule


\begin{example}\label{e3} {\bf $(3^5, 4)$}-SEM on the double torus. These example are constructed by lifting the SEMs $K_i$ to its double cover $T_i$, $i.e.$, the double torus. Consider a map defined by $\phi\{0, 12\} = a$; $\phi\{1, 18\} = b$; $\phi\{2, 20\} = c$; $\phi\{3, 21\} = d$; $\phi\{4, 19\} = e$; $\phi\{5, 13\} = f$; $\phi\{6, 23\} = g$; $\phi\{7, 17\} = h$; $\phi\{8, 18\} = i$; $\phi\{9, 15\} = j$; $\phi\{10, 14\} = k$; $\phi\{11, 22\} = l$ and the map $\psi\colon \{a \mapsto 0, b\mapsto 1 , \ldots,l \mapsto 11\}$. We see that $\psi\circ\phi \colon T_{i} \longrightarrow K_{i}$, $i = 1, 2, 3$ is a covering. Clearly it is a two fold orientable covering\,:


 {\begin{flushleft}{$\bf T_1$} := \{\rm{[0, 1, 2],  [0, 1, 7], [0, 6, 7], [0, 6, 13], [0, 4, 13], [1, 2, 8], [1, 5, 8], [1, 5, 10], [2, 7, 8], [2, 6, 7], [2, 3, 6],[3, 6, 9], [3, 9, 10], [3, 10, 22], [3, 4, 22], [4, 22, 15], [4, 15, 14], [4, 14, 13], [5, 10, 19], [5, 19, 12], [5, 12, 23],  [7, 8, 22], [8, 15, 22], [9, 10, 19], [9, 11, 16], [9, 11, 19], [11, 16, 17], [11, 19, 21], [11, 14, 21], [12, 17, 23], [12, 17, 18],  [12, 18, 20], [13, 14, 18], [13, 16, 18], [14, 15, 21], [15, 21, 23], [16, 17, 20], [16, 18, 20],  [17, 20, 23],  [20, 21, 23]}\} $\bigcup$ \{\rm{[0, 2, 3, 4], [1, 7, 22, 10], [5, 8, 15, 23], [6, 9, 16, 13], [11, 14, 18, 17],  [12, 19, 21, 20]}\}
 \end{flushleft}} 


{\begin{flushleft}{$\bf T_2$} := \{\rm{[0, 2, 13],  [0, 13, 7], [0, 6, 7], [0, 5, 6], [0, 4, 5], [1, 8, 9], [1, 9, 14], [1, 12, 14], [1, 12, 19],
[1, 11, 19], [2, 3, 8],  [2, 6, 8], [2, 6, 21], [2, 13, 21], [3, 8, 9], [3, 9, 22], [3, 22, 23], [3, 4, 23], [4, 5, 10], [4, 7, 10],
[4, 7, 23], [5, 6, 8], [5, 10, 11], [7, 13, 23], [9, 14, 18], [10, 11, 15], [10, 15, 21], [11, 15, 16], [11, 16, 19], [12, 18, 19],
[12, 17, 18], [12, 16, 17], [13, 20, 21], [14, 15, 20], [14, 18, 20], [15, 20, 21], [16, 17, 22], [16, 19, 22], [17, 18, 20],
[17, 22, 23]}\} $\bigcup$ \{\rm {[0, 2, 3, 4], [1, 8, 5, 11], [6, 7, 10, 21], [9, 18, 19, 22], [12, 14, 15, 16], [13, 20, 17, 23]} \}\end{flushleft}} 


{\begin{flushleft} {$\bf T_3$} := \{\rm{[0, 1, 2],  [0, 1, 7], [0, 6, 7], [0, 5, 6], [0, 4, 5], [1, 2, 9], [1, 9, 11], [1, 7, 8], [2, 3, 8], [2, 6, 8],
[2, 6, 9],  [3, 7, 8], [3, 7, 10], [3, 10, 23], [3, 4, 23], [4, 5, 22], [4, 21, 22], [4, 21, 23], [5, 6, 8], [5, 11, 22], [9, 11, 16],
[9, 10, 16], [10, 11, 17], [10, 16, 17], [11, 15, 16], [11, 15, 22], [12, 16, 17], [12, 17, 18], [12, 18, 19], [12, 13, 19], [12, 13, 14], [13, 14, 21], [13, 21, 23], [13, 19, 20], [14, 18, 20], [14, 18, 21], [14, 15, 20], [15, 19, 22], [15, 19, 20], [17, 18, 20]}\} $\bigcup$ \{\rm{[0, 2, 3, 4], [1, 8, 5, 11], [6, 7, 10, 9], [12, 14, 15, 16], [13, 20, 17, 23], [18, 19, 22, 21]}\}
\end{flushleft}}

\end{example}


\begin{example}\label{eN} Following is the example of a SEM of type $(3^5, 4)$ on a non - orientable surface of Euler characteristics $-2$\,:

{\begin{flushleft} {$\bf N$} := \{\rm{[0, 1, 2], [0, 1, 18], [0, 18, 14], [0, 14, 15], [0, 15, 4],[1, 2, 8],[1, 8, 5], [1, 5, 10], [2, 3, 6], [2, 6, 7], [2, 7, 8], [3, 6, 13], [3, 10, 13], [3, 10, 11],[3, 4, 11], [4, 11, 9], [4, 9, 16], [4, 15, 16], [5, 1, 8], [5, 10, 17], [5, 17, 23], [5, 6, 23], [6, 7, 23], [7, 8, 12], [7, 22, 23], [8, 12, 13], [9, 11, 19], [9, 16, 21], [9, 14, 21], [10, 13, 17], [12, 13, 17], [12, 17, 21], [12, 21, 16], [14, 18, 20], [14, 20, 21], [15, 16, 22], [15, 19, 22], [18, 19, 20], [18, 11, 19], [19, 20, 22], [20, 22, 23], [0, 2, 3, 4], [1, 10, 11, 18], [5, 6, 13, 8], [7, 12, 16, 22], [9, 14, 15, 19], [17, 21, 20, 23]}\} \end{flushleft}}

\end{example}

\begin{table}
\centering
\caption {\bf Table representing cylinder additions to $K_{ij}$s}\label{Table1}
\smallskip
{\small{
\begin{tabular}{|l|l|l|l|l|l|l|}
   \hline
  $Maps $ & $\#1$ & $\#5$ & $\#6$ &$\chi$ & Or.& Handle Type\\
  \hline
  $K_{11}(1, 1)$ & 0 & 10 & 4 & -8 & NO & $C_{44}$([0, 2, 3, 4],[18, 21, 19, 20]), $C_{44}$([1, 7, 11, 10],\\&&&&&& [16, 23, 13, 12]),
  $C_{44}$([5, 6, 9, 8], [14, 15, 22, 17])\\
  \hline
  $K_{11}(2, 1)$ & 26 & 36 & 2 & -10 & NO & $C_{33}$([0, 1, 2], [12, 13, 14]),
  $C_{33}$([ 3, 6, 9 ], [15, 18, 21])\\
  &&&&&&$C_{33}$([ 4, 5, 10 ], [16, 17, 22]), $C_{33}$([ 7, 8, 11 ], [19, 20, 23])\\
  \hline
   $K_{22}(1, 1)$ & 0 & 12 & 0 & -8 & NO & $C_{44}$([0, 2, 3, 4], [12, 14, 15, 16]), $C_{44}$([1, 8, 5, 11],\\&&&&&&  [13, 20, 17, 23]), $C_{44}$([6, 9, 10, 7], [21, 18, 19, 22])\\
   \hline
   $K_{22}(2, 1)$ & 16 & 32 & 4 & -10 & NO & $C_{33}$([0, 1, 2], [12, 13, 14]), $C_{33}$([ 3, 9, 10 ], [15, 21, 22]),\\
   &&&&&&$C_{33}$([ 4, 7, 11 ], [16, 19, 23]), $C_{33}$([ 5, 6, 8 ], [17, 18, 20])\\
   \hline
  $K_{33}(1, 1)$ & 0 & 16 & 0 & -8 & NO & $C_{44}$([0, 2, 3, 4], [12, 14, 15, 16]), $C_{44}$([6, 7, 10, 9],\\ &&&&&& [19, 22, 21, 18]),
  $C_{44}$([1, 11, 5, 8],  [13, 23, 17, 20])\\
  \hline
   $K_{33}(2, 1)$ & 18 & 34 & 6 & -10 & NO & $C_{33}$([0, 1, 2], [12, 13, 14]), $C_{33}$([ 3, 7, 10 ], [15, 19, 22]),\\
   &&&&&&$C_{33}$([ 4, 9, 11 ], [16, 21, 23]),$C_{33}$([ 5, 6, 8 ], [17, 18, 20])\\
   \hline
   $K_{12}(1, 1)$ & 1 & 10 & 2 & -8 & NO & $C_{33}$([0, 2, 3, 4],  [13, 14, 15, 23]), $C_{33}$([1, 7, 11, 10],\\ &&&&&& [19, 16, 22, 12]), $C_{33}$([5, 8, 9, 6],  [17, 18, 21, 20])\\
   \hline
  $K_{12}(2, 1)$ & 19 & 34 & 3 & -10 & NO & $C_{33}$([0, 1, 2 ], [ 12, 13, 14]), $C_{33}$([3, 6, 9 ], [ 15, 21, 22 ]),\\
  &&&&&&$C_{33}$([ 4, 5, 10 ], [ 16, 19, 23 ]), $C_{33}$([ 7, 8, 11 ], [ 17, 18, 20 ])\\
  \hline
   $K_{13}(1, 1)$ & 1 & 15 & 2 & -8 & NO & $C_{44}$([0, 2, 3, 4],  [12, 14, 15, 16]), $C_{44}$([1, 7, 11, 10],\\ &&&&&&  [13, 23, 17, 20]), $C_{44}$([5, 8, 9, 6],  [18, 19, 22, 21])\\
   \hline
  $K_{13}(2, 1)$& 21 & 35 & 4 & -10 & NO & $C_{33}$([0, 1, 2], [ 12, 13, 14]), $C_{33}$([3, 6, 9 ], [ 15, 19, 22 ]),\\
  &&&&&&$C_{33}$([ 4, 5, 10 ], [ 16, 21, 23 ]), $C_{33}$([ 7, 8, 11 ], [ 17, 18, 20 ])\\
  \hline
   $K_{23}(1, 1)$ & 0 & 15 & 0 & -8 & NO & $C_{44}$([0, 2, 3, 4], [12, 14, 15, 16]), $C_{44}$([1, 8, 5, 11],\\ &&&&&& [13, 23, 17, 20]), $C_{44}$([6, 9, 10, 7], [18, 21, 22, 19])\\
   \hline
   $K_{23}(2, 1)$ & 15 & 33 & 5 & -10 & NO & $C_{44}$([0, 1, 2], [12, 13, 14]), $C_{33}$([3, 9, 10 ], [ 15, 19, 22 ]),\\
   &&&&&&$C_{33}$([4, 7, 11 ], [ 16, 21, 23 ]), $C_{33}$([5, 6, 8 ], [17, 18, 20 ]).\\
  \hline
\end{tabular}}}

\smallskip

  \centering
  \caption{\bf Table representing Cylinder Additions to double covers}\label{Table2}
  \smallskip
{\small{
\begin{tabular}{|l|l|l|l|l|l|l|}
  \hline
  ${ Maps}$ & $E_1$ & $E_5$ & $E_6$ & $ \chi$ & Or. & Handle Type\\
  \hline
   $T_1(1, 1)$& 4 & 18 & - & -8 & NO & $C_{44}$([24, 2, 3, 4], [11, 14, 18, 17]),$C_{44}$([12, 19, 21, 20],\\ &&&&&& [1, 7, 22, 10]), $C_{44}$([5, 8, 15, 23], [6, 9, 16, 13])\\
   \hline
  $T_1(1, 2)$ & 14& 24 & 2 & -8 & NO & $C_{44}$([24, 2, 3, 4], [11, 14, 18, 17]), $C_{44}$([12, 19, 21, 20],\\ &&&&&& [1, 7, 22, 10]), $C_{44}$([5, 8, 15, 23], [6, 9, 16, 13]), \\
      \hline
  $T_2(1, 1)$ & 4 & 26 & - & -8 & O & $C_{44}$([6, 7, 10, 21], [9, 18, 19, 22]), $C_{44}$([24, 2, 3, 4],\\ &&&&&& [12, 14, 15, 16]), $C_{44}$([13, 20, 17, 23],[1, 8, 5, 11]) \\
  \hline
  $T_2(2, 1)$ & 11 & 21 & 4 & -10 & NO & $C_{33}$([ 24, 2, 13 ], [ 11, 16, 19 ]), $C_{33}$([ 1, 9, 14 ], [ 4, 7, 10 ]),\\
    &&&&&&$C_{33}$([ 3, 22, 23 ], [ 15, 20, 21 ]), $C_{33}$([ 5, 6, 8 ], [ 12, 17, 18 ]) \\
  \hline
  $T_3(1, 1)$ & 0 & 24 & - & -8 & O & $C_{44}$([24, 2, 3, 4], [12, 14, 15, 16]), $C_{44}$([1, 8, 5, 11],\\ &&&&&& [13, 20, 17, 23]), $C_{44}$([18, 19, 22, 21], [6, 7, 10, 9])\\
  \hline
  $T_3(2, 1)$ & 9 & 17 & 6 & -10 & O & $C_{33}$([ 24, 1, 2 ], [ 12, 13, 14 ]), $C_{33}$([ 3, 7, 10 ], [ 15, 19, 22 ]),\\
   &&&&&&$C_{33}$([ 4, 21, 23 ], [ 9, 11, 16 ]), $C_{33}$([ 5, 6, 8 ], [ 17, 18, 20 ])\\
  \hline
  $T_3(2, 2)$ & 20 & 23 & 10 & -10 & NO & $C_{33}$([ 24, 1, 2 ],[ 17, 18, 20 ]), $C_{33}$([ 3, 7, 10 ], [ 12, 13, 14 ]),\\
    &&&&&&$C_{33}$([ 4, 21, 23 ], [ 9, 11, 16 ]), $C_{33}$([ 5, 6, 8 ], [ 15, 19, 22 ])\\
  \hline
  $T_3(2, 3)$ & 11 & 17 & 6 & -10 & NO & $C_{33}$([ 24, 1, 2 ], [ 17, 18, 20 ]), $C_{33}$([ 3, 7, 10 ], [ 15, 19, 22 ]),\\
   &&&&&&$C_{33}$([ 4, 21, 23 ], [ 9, 11, 16 ]), $C_{33}$([ 5, 6, 8 ], [ 12, 13, 14 ]) \\
  \hline
  $T_3(2, 4)$ & 21 & 18 & 12 & -10 & O & $C_{33}$([ 24, 1, 2 ], [ 15, 19, 22 ]), $C_{33}$([ 3, 7, 10 ],[ 12, 13, 14 ]),\\
    &&&&&&$C_{33}$([ 4, 21, 23 ], [ 9, 11, 16 ]), $C_{33}$([ 5, 6, 8 ], [ 17, 18, 20 ]\\
  \hline
\end{tabular}}}
\end{table}

\begin{example}\label{e4} In the Table \ref{Table1} some SEMs of type {\bf $(3^5, 4^2)$} and $(3^7, 4)$ on the surfaces of Euler characteristic $-8$ and $-10$ are presented. These examples are obtained from $K_i$s  by cylinder addition techniques. The last column in the tables gives the faces where cylinders are added. The notation $K_{ij}(k, l)$ denotes $l^{th}$ example in the set of objects obtained by adding quadrangle for $k = 1$ and triangle for $k = 2$ between $K_i$ and $K_j$.
\end{example}

\begin{lem}\label{lemma1} The maps defined in the Example \ref{e4} above are all non - isomorphic.
\end{lem}
\begin{proof} Consider the enumeration of edges in $G_1$, $G_5$ and $G_6$ of the SEMs in this example presented in tabular form above. From this it is immediate that the SEMs of this example are all non-isomorphic.
\end{proof} $\hfill \Box$
\begin{example} \label{e5} Table \ref{Table2} presents some SEM which are obtained by adding cylinders in the double covers $T_1$, $T_2$ and $T_3$. The notation $T_i(j, k)$ denotes the $k^{th}$ object obtained from the double torus $T_i$ of example \ref{e3}, by adding the cylinders of type $C_{44}$ for $j = 1$ and of type $C_{33}$ for $j = 2$. The last column in the tables gives the faces where cylinders are added.
\end{example}

\begin{lem}\label{lemma2} The maps defined in Example \ref{e5} above are mutually non isomorphic.
\end{lem}
\begin{proof} Consider the enumeration of edges in $G_1$, $G_5$ and $G_6$ of the SEMs in this example presented in tabular form above. From this it is immediate that the SEMs of this example are all non-isomorphic.
\end{proof} $\hfill \Box$

\section{Proofs}

In this section we present proof of the results given in introduction section.
\smallskip

\noindent {\sc Proof of Theorem \ref{t2}}:\, The proof follows by examples \ref{e3} and \ref{eN}. It is obvious that $N$ is not isomorphic to any of $T_1$, $T_2$ and $T_3$. Now, $EG(G_5(T_1)) = \{ [ 1, 7 ], [ 2, 22 ], [ 2, 24 ], [ 3, 7 ], [ 11, 20 ], [ 12, 20 ] \}$ and $EG(G_6(T_1)) = \emptyset$.  $EG(G_5(T_2)) = \{ [ 2, 7 ], [ 4, 6 ], [ 14, 19 ], [ 16, 18 ] \}$ and $EG(G_6(T_2)) = \emptyset$.  $EG(G_6$ $(T_3))$= $\{$ [ 1, 6 ], [ 2, 7 ], [ 8, 24 ], [ 12, 20 ], [ 13, 18 ], [ 14, 19 ] $\}$. Hence $T_1 \not\cong T_2$, $T_1 \not\cong T_3$ and $T_2 \not\cong T_3$. From here it is also evident that $T_1$, $T_2$ and $T_3$ are not vertex transitive. Also, since $ EG(G_4(N)) = \{[ 1, 3 ], [ 5, 13 ], [ 6, 8 ], [ 9, 15 ], [ 11, 24 ], [ 12, 23 ], [ 14, 19 ], [ 21, 22 ]\}$ it follows that $N$ is also not vertex transitive. $\hfill\Box$


\noindent {\sc Proof of Theorem \ref{t3}:} The result follows from Lemma \ref{lemma1} and Lemma \ref{lemma2}. \hfill$\Box$


\noindent {\sc Proof of Theorem \ref{t4}:} The result follows from Lemma \ref{lemma1} and Lemma \ref{lemma2}. \hfill$\Box$


\noindent {\sc Proof of Theorem \ref{t1}}\, Let $K$ be a SEM of type $(3^5, 4)$ on the surface of Euler characteristic $-1$. Let $V = V(K) = \{0, 1, 2, \ldots, 10, 11\}$ denote the set of vertices of $K$. The proof of the theorem is by exhaustive search for all $K$. In what follows, the notation ${\rm lk}(i) = C_7([i_1, i_2, i_3], i_4, i_5, i_6, i_7)$ for link of $i$ will mean that $[i, i_1, i_2, i_3]$ forms a quadrangular facet and $[i, i_3, i_4]$, $[i, i_4, i_5]$, $[i, i_5, i_6]$, $[i, i_6, i_7]$, $[i, i_7, i_1]$ form triangular facets. We may assume without loss of generality that ${\rm lk}(0) = C_7([2 ,3, 4], 5, 6, 7, 1)$. Then ${\rm lk}(2) = ([3,4,0],1,c,b,a)$. It is easy to see that $(a, b, c) \in A\cup B\cup C$, where

{\begin{flushleft} ${\bf A} = \{$(6, 7, 8),  (8, 6, 9), (8, 7, 6), (8, 9, 6), (8, 9, 10)$\}$, ${\bf B} = \{$(5, 6, 7), (5, 6, 8), (5, 7, 6), (5, 7, 8),(6, 5, 7), (6, 5, 8), (5, 8, 6), (5, 8, 9), (6, 7, 5),  (6, 8, 5),
(6, 8, 9), (7, 8, 5), (7, 8, 6), $\}$ and ${\bf C} = \{$ (7, 5, 6), (7, 5, 8), (7, 6, 8), (7, 8, 9), (8, 5, 6),  (8, 5, 9), (8, 6, 5), (8, 7, 5), (8, 7, 9), (8, 9, 5)$\}$ \end{flushleft}}


\begin{claim}\label{cl1}The values $(a, b, c) \in C$ are isomorphic to some $(a, b, c) \in A\cup B$.
\end{claim}
{\begin{flushleft} We have $(0, 2)(3, 4)(5, 8)(6, 7, 9) \colon (8, 9, 6) \cong (8, 7, 9)$~~$(5, 7, 9, 6, 8) \colon (8, 9, 5) \cong (5, 6, 7)$~~$(5, 7, 6, 8) \colon (8, 7, 5) \cong (5, 6, 7)$~~$(5, 8, 6, 9, 7) \colon (7, 8, 9) \cong (5, 6, 7)$~~$(5, 6, 7, 8) \colon (8, 5, 6) \cong (5, 6, 7)$~~$(5, 6, 8, 7) \colon (7, 5, 8) \cong (5, 6, 7)$~~$(5, 7, 8) \colon (8, 6, 5) \cong (5, 6, 7)$~~$(5, 6, 8)(7, 9) \colon (8, 5, 9) \cong (5, 6, 7)$~~$(0, 2)(3, 4)(7, 5, 6) \colon (7, 5, 6) \cong (6, 7, 5)$~~$(0, 2)(3, 4)(7, 5, 8) \colon (7, 6, 8) \cong (8, 6, 5)$. This proves the Claim \ref{cl1}.

Thus we have $S = \{(6, 7, 8)$, $(8, 6, 9)$, $(8, 7, 6)$, $(8, 9, 6)$, $(8, 9, 10)\}$ $\cup\{(5, 6, 8)$, $(5, 7, 6)$, $(5, 7, 8)$, $(6, 5, 7)$, $(6, 5, 8)$, $(5, 8, 6)$, $(5, 8, 9)$, $(6, 7, 5)$, $(6, 8, 5)$, $(6, 8, 9)$, $(7, 8, 5)$, $(7, 8, 6)\}$.
\end{flushleft}}

\begin{claim}\label{cl2} There are no semi-equivelar maps(SEM) on a surface of Euler characteristic -1 for the values $(a, b, c) \in B$.
\end{claim}

\smallskip

\noindent By considering ${\rm lk}(1)$ we see immediately that if $(a, b, c) = (5, 6, 7)$ then $C_3(2, 0, 7) \subseteq {\rm lk}(1)$. When $(a, b, c) = (5, 6, 8)$ this implies 3 or 4 appears in the squares containing 5, $i.e.$ they appear in two square. This is not possible. When $(a, b, c) = (5, 7, 6)$ then ${\rm lk}(2) = C_7([3, 4, 0], 1, 6, 7, 5)$ this implies ${\rm lk}(5) = ([7, a, 6], 0, 4, 3, 2)$ then $C_4(3, 5, 0, 4)\subseteq{\rm lk}(2)$. Which is a contradiction. When $(a, b, c) = (5, 7, 8)$ as in previous case we observe that this case is not possible. When $(a, b, c) = (6, 5, 7)$ then ${\rm lk}(2) = C_7([3, 4, 0], 1, 7, 5, 6)$. Considering links of $6, 7, $ and $5$ successively we get a contradiction as in case of $(a, b, c) = (5, 6, 8)$. When $(a, b, c) = (6, 5, 8)$ then ${\rm lk}(6) = C_7([7, 8, 9], 3, 2, 5, 0)$. Considering link of $5$ we get a contradiction as in the case $(a, b, c) = (5, 6, 8)$. When $(a, b, c) = (5, 8, 6)$ then considering links of $2$ and $5$ successively we see that $\deg(4) \leq 4$. Similarly, when $(a, b, c) = (5, 8, 9)$ as in previous case this case also leads to a contradiction. When $(a, b, c) = (6, 7, 5)$ then considering links of 5 and 7 successively we get a contradiction as in case $(a, b, c) = (5, 6, 8)$. When $(a, b, c) = (6, 8, 5)$ then considering link 2 we get $\deg(5) \geq 8$.

\smallskip

When $(a, b, c) = (6, 8, 9)$ then ${\rm lk}(2) = C_7([3, 4, 0], 1, 9, 8, 6)$ this implies ${\rm lk}(6) = C_7([8, 10, 7], 0, 5$, $3, 2)$ or ${\rm lk}(6) = C_7([8, 10, 5], 0, 7, 3, 2)$. In the first case, successively considering links of $3$, $5$, $7$, $4$, $11$, $9$ and $1$ we see that this case is not possible. In the second case, ${\rm lk}(7) = C_7([9, 11, 1], 0, 6, 3, x)$ or ${\rm lk}(7) = C_7([11, 9, 1], 0, 6, 3, x)$, for some $x \in V$. In the first case, $19$ is both a non-edge and an edge in $M$ and in second case, $C_5(9, 2, 0, 7, 11) \subseteq {\rm lk}(1)$. A contradiction. When $(a, b, c) = (7, 8, 5)$ then considering links of 2 and 5 we get a contradiction as in case $(a, b, c) = (5, 6, 8)$. When $(a, b, c) = (7, 8, 6)$ then link 6 has more than 7 vertices.

\smallskip

When $(a, b, c) = (7, 8, 9)$ then ${\rm lk}(2) = C_7([3, 4, 0], 1, 9, 8, 7)$ this implies ${\rm lk}(7) = C_7([8, 10, 1], 0, 6$, $3, 2)$ or ${\rm lk}(7) = C_7([8, 10, 6], 0, 1, 3, 2)$. When ${\rm lk}(7) = C_7([8, 10, 1], 0, 6, 3, 2)$ we get ${\rm lk}(1) = C_7([10, 8, 7], 0, 2, 9, d)$ where $d \in \{4, 5, 11\}$. When $d = 4$ then ${\rm lk}(4) = C_7([3, 2, 0], 5, 9, 1, 10)$ or ${\rm lk}(4) = C_7([3, 2, 0], 5, 10, 1, 9)$. First case implies ${\rm lk}(9) = C_7([11, 6, 5], 4, 1, 2, 8)$ then $C_5(4, 0, 6, 11, 9)$ $\subseteq {\rm lk}(5)$. In second case considering ${\rm lk}(9)$ we see that 8 or 3 appear in two square. When $d = 5$ then ${\rm lk}(1) = C_7([10, 8, 7], 0, 2, 9, 5)$, this implies ${\rm lk}(5) = C_7([9, 11, 6]$, $0, 4, 10, 1)$, ${\rm lk}(6) = C_7([11, 9, 5], 0, 7, 3, e)$ where $e \in\{4, 8\}$. If $e = 4$ then ${\rm lk}(4) = C_7([0, 2, 3], 6, 11, 10, 5)$ this implies $C_5(4, 6, 7, 2, 0)\subseteq{\rm lk}(3)$. If $e = 8$, ${\rm lk}(8)$ has more than seven vertices. When $d = 11$ we get ${\rm lk}(3) = C_7([4, 0, 2], 7, 6, x, y)$ there $(x, y)\in\{(8, 9), (9, 8)$, $(8, 10), (10, 8), (10, 11), (11, 10), (9, 11), (11, 9)\}$. If $(x, y)\in\{(8, 9), (9, 11), (11, 9)\}$ then we do not get three quadrangles in the map. So, $(x, y)\in\{(9, 8), (8, 10), (10, 8), (10, 11), (11, 10)\}$.

If $(x, y)= (9, 8)$ then ${\rm lk}(3) = C_7([4, 0, 2], 7, 6, 9, 8)$. So, ${\rm lk}(9) = C_7([11, 5, 6], 3, 8, 2, 1)$. This implies $C_6(5, 0, 7, 3, 9, 11) \subseteq{\rm lk}(6)$. Which is not possible. If $(x, y)= (8, 10)$ then ${\rm lk}(3) = C_7([4, 0, 2], 7, 6, 8, 10)$, ${\rm lk}(8) = C_7([7, 1, 10], 3, 6, 9, 2)$, ${\rm lk}(6) = C_7([9, 11, 5], 0, 7, 3, 8)$. Then ${\rm lk}(9)$ has a 6-cycle. If $(x, y)= (10, 8)$ then ${\rm lk}(3) = C_7([4, 0, 2], 7, 6, 10, 8)$ and ${\rm lk}(8) = C_7([7, 1, 10], 3, 4, 9, 2)$. This implies ${\rm lk}(6) = C_7([5, 9, 11], 10, 3, 7, 0)$. This implies ${\rm lk}(10)$ has a 6-cycle. This is not allowed. If $(x, y)= (10, 11)$ then ${\rm lk}(3) = C_7([4, 0, 2], 7, 6, 10, 11)$ and ${\rm lk}(10) =C_7([8, 7, 1], 11, 3, 6, f)$. Possible values of $f = 5$ or 9. If $f = 5$, then $C_5(5, 10, 3, 7, 0)\subseteq{\rm lk}(6)$. If $f = 9$ then $C_5(9, 10, 1, 7, 2)\subseteq{\rm lk}(8)$. A contradiction. If $(x, y)= (11, 10)$ then ${\rm lk}(3) = C_7([4, 0, 2], 7, 6, 11, 10)$. This implies ${\rm lk}(10) = C_7([8, 7, 1], 11, 3, 4, g)$, for some $g \in V$. It is easy to see that $g = 5$ or 9. When $g = 5$, we get $C_5(5, 10, 3, 2, 0)\subseteq{\rm lk}(4)$. When $g = 9$ then we get $C_5(9, 10, 1, 7, 2)\subseteq{\rm lk}(8)$. This is not possible. So, $d \neq 11$. Hence ${\rm lk}(7) = C_7([8, 10, 6], 0, 1, 3, 2)$. In this case, we get ${\rm lk}(1) = C_7([5, 11, 9], 2, 0, 7, 3)$ or ${\rm lk}(1) = C_7([11, 5, 9], 2, 0, 7, 3)$.

\smallskip

When ${\rm lk}(1) = C_7([5, 11, 9], 2, 0, 7, 3)$ we get ${\rm lk}(5) = C_7([1, 9, 11], 4, 0, 6, 3)$, ${\rm lk}(6) = C_7([7, 8, 10]$, $4$, $3$, $5$, $0)$, ${\rm lk}(4) = C_7([3, 2, 0], 5, 11, 10, 6)$, ${\rm lk}(3) = C_7([4, 0, 2], 7, 1, 5, 6)$. So, ${\rm lk}(11) = C_7([9, 1, 5]$, $4$, $10$, $w$, $z)$. But then $\{0, 2, 3, 6, 7\}$ does not belongs to ${\rm lk}(11)$. That is $\deg(11) \leq 6$. When ${\rm lk}(1) = C_7([11, 5, 9], 2, 0, 7, 3)$ we get ${\rm lk}(3) = C_7([4, 0, 2], 7, 1, 11, 10)$. This implies ${\rm lk}(4) = C_7(10, [3, 2, 0]$, $5, p, q)$, for some $p, q \in V$. It is easy to see that $(p, q)\in\{(8, 9)$, $(8, 11)$, $(9, 8)$, $(11, 8)\}$. If $(p, q) = (8, 11)$, then ${\rm lk}(8)$ has more than seven vertices. If $(p, q) = (9, 8)$ then $C_(5, 4, 8, 2, 1, 11)\subseteq {\rm lk}(9)$. When $(p, q) = (11, 8)$ then ${\rm lk}(4) = C_7([0, 2, 3], 10, 8, 11, 5)$. So, ${\rm lk}(8) = C_7([10, 6, 7], 2, 9, 11, 4)$. This implies $28$ and $811$ are edges in ${\rm lk}(9)$ which is not possible. When $(p, q) = (8, 9)$ then ${\rm lk}(4) = C_7([0, 2, 3], 10, 9, 8, 5)$, ${\rm lk}(8) = C_7([10, 6, 7], 2, 9, 4, 5)$, ${\rm lk}(9) = C_7([1, 11, 5], 10, 4, 8, 2)$, ${\rm lk}(10) = C_7([8, 7, 6], 3, 4, 9, 5)$ and ${\rm lk}(5) = C_7([11, 1, 9], 10, 8, 4, 0)$. This implies ${\rm lk}(6) = C_7([7, 8, 10], 3, 11$, $5, 0)$. But then, $C_3(3, 11, 10) \in {\rm lk}(3)$ which is a contradiction. This proves the Claim \ref{cl2}.


Thus we may assume $(a, b, c) \in A$. In other words $(a, b, c) \in \{(6, 7, 8)$,  $(8, 6, 9)$, $(8, 7, 6)$, $(8, 9, 6)$, $(8, 9, 10)\}$.


\smallskip

\noindent{\bf Case 1:} When $(a, b, c) = (6, 7, 8)$,  ${\rm lk}(2) = C_7([3, 4, 0], 1, 8, 7, 6)$. So, ${\rm lk}(6) = C_7([8, 1, 5], 0, 7, 2, 3)$, ${\rm lk}(6) = C_7([8, 9, 5], 0, 7, 2, 3)$, ${\rm lk}(6) = C_7([1, 8, 5], 0, 7, 2, 3)$ or ${\rm lk}(6) = C_7([9, 8, 5], 0, 7, 2, 3)$. When ${\rm lk}(6) = C_7([8, 1, 5], 0, 7, 2, 3)$, considering ${\rm lk}(7)$, we get 8 in two quadrangles.  When ${\rm lk}(6) = C_7([8, 9, 5], 0, 7, 2, 3)$, considering links of $8$ and $7$ we get either 3 or 9 in two quadrangles. This is not possible. Similarly, when ${\rm lk}(6) = C_7([1, 8, 5], 0, 7, 2, 3)$ then either 1 or 8 will be in two squares. Thus ${\rm lk}(6) = C_7([9, 8, 5], 0, 7, 2, 3)$. This implies ${\rm lk}(7) = C_7([11, 10, 1], 0, 6, 2, 8)$. This implies ${\rm lk}(8) = C_7([5, 6, 9], 1, 2, 7, 11)$ or ${\rm lk}(8) = C_7([5, 6, 9], 11, 7, 2, 1)$. If ${\rm lk}(8) = ([5, 6, 9], 1, 2, 7, 11)$ then considering links of 1 and 5 we get $C_5(10, 5, 8, 7, 1)\subseteq{\rm lk}(11)$. This is not possible. So, ${\rm lk}(8) = C_7([5, 6, 9], 11, 7, 2, 1)$ now completing successively we get ${\rm lk}(1) = C_7([7, 11, 10], 5$, $8, 2, 0)$, ${\rm lk}(5) = C_7([6, 9, 8], 1, 10, 4, 0)$, ${\rm lk}(10) = C_7([1, 7, 11], 3, 9, 4, 5)$, ${\rm lk}(9) = C_7([6, 5, 8]$, $11, 4, 10$, $3)$, ${\rm lk}(3) = C_7([2, 0, 4], 11, 10, 9, 6)$, ${\rm lk}(4) = C_7([3, 2, 0], 5, 10, 9, 11)$, ${\rm lk}(11) = C_7([10, 1, 7], 8, 9, 4, 3)$. This is the object \noindent{\bf $K_1$\,} of Example \ref{e1} with vertex 10 replaced by $u$ and vertex 11 replaced by $v$.

\smallskip

\noindent{\bf Case 2:} When $(a, b, c) = (8, 6, 9)$ then ${\rm lk}(2) = C_7([3, 4, 0], 1, 9, 6, 8)$. This implies ${\rm lk}(6) = C_7(0, 7, 9, 2, [8, d, 5])$, ${\rm lk}(6) = C_7(0, 7, 8, 2, [9, d, 5])$,  ${\rm lk}(6) = C_7(0, 5, 9, 2, [8, d, 7])$ or ${\rm lk}(6) = C_7(0, 5, 8$, $2, [9, d, 7])$, for some $d \in V$. If ${\rm lk}(6) = C_7(0, 7, 8, 2, [9, d, 5])$, then it is easy to see that $d = 1$ or $10$. If $d = 1$, then $C_4(1, 2, 6, 5)\subseteq{\rm lk}(9)$. If $d = 10$ then for some $x, y, z \in V$ we have ${\rm lk}(7) = C_7(0, 6, 8, x, [y, z, 1])$ or ${\rm lk}(7) = C_7(6, 0, 1, x, [y, z, 8])$. But in both cases ${\rm lk}(7)$ can not be completed. If ${\rm lk}(6) = C_7(0, 7, 9, 2, [8, d, 5])$, then $d = 1$ or $10$ If $d = 1$, we get ${\rm lk}(7) = C_7([11, 10, 9], 6, 0, 1, e)$ or ${\rm lk}(7) = C_7([10, 11, 9], 6, 0, 1, e)$, for some $e \in V$.
In the first case, $i.e.$ when ${\rm lk}(7) = C_7([11, 10, 9], 6, 0, 1, e)$ we see that $e \in \{3, 4, 5, 8\}$. If $e = 3$ or 4 then ${\rm lk}(1)$
has more than seven vertices. When $e = 5$ then we get ${\rm lk}(1) = C_7([5, 6, 8], 9, 2, 0, 7)$, ${\rm lk}(8) = C_7([6, 5, 1], 9, 10, 3, 2)$,
 ${\rm lk}(9) = C_7([7, 11, 10], 8, 1, 2, 6)$ and ${\rm lk}(5) = C_7([6, 8, 1], 7, 11, 4, 0)$. But then ${\rm lk}(11)$ can not be completed since
 $\{0, 1, 2, 6, 8\}$ does not belong to ${\rm lk}(11)$. So, $e = 8$. Then ${\rm lk}(7) = C_7([11, 10, 9], 6, 0, 1, 8)$. This implies ${\rm lk}(9) = C_7([7, 11, 10], 5, 1, 2, 6)$,  ${\rm lk}(5) = C_7([6, 8, 1], 9, 10, 4, 0)$. Then considering ${\rm lk}(10)$ we see that $C_5(3, 10, 5, 0, 2)\subseteq{\rm lk}(4)$. This is not possible. When ${\rm lk}(7) = C_7([10, 11, 9], 6, 0, 1, e)$,  we see that $e\in\{3, 4, 5, 8\}$. If $e = 3$ or 4 then ${\rm lk}(1)$
has more than seven vertices. If $e = 5$, then ${\rm lk}(1) = C_7([5, 6, 8], 9, 2, 0, 7)$. This implies ${\rm lk}(9) = C_7([7, 10, 11], 8, 1, 2, 6)$
${\rm lk}(8) = C_7([6, 5, 1], 9, 11, 3, 2)$ and ${\rm lk}(6) = C_7([5, 1, 8], 2, 9, 7, 0)$. Then $\{0, 1, 2, 4, 5\} \not\in {\rm lk}(11)$. This is not possible. If $e = 8$ then ${\rm lk}(7) = C_7([10, 11, 9], 6, 0, 1, 8)$, ${\rm lk}(6) = C_7([8, 1, 5], 0, 7, 9, 2)$, ${\rm lk}(8) = C_7([1, 5, 6], 2, 3, 10, 7)$, ${\rm lk}(1) = C_7([5, 6, 8], 7,  0, 2, 9)$, ${\rm lk}(5) = C_7([6, 8, 1], 9, 11, 4, 0)$ and ${\rm lk}(9) = C_7([7, 10, 11], 5, 1, 2, 6)$. Then $\{0, 1, 2, 6, 8\} \not\in {\rm lk}(11)$. This is not possible. Thus, $d \neq 1$. Hence $d = 10$. Then ${\rm lk}(6) = C_7([5, 10, 8], 2, 9, 7, 0)$. This implies that for some $x, y, z \in V$ we have ${\rm lk}(7) = C_7(0, 6, 9, x, [y, z, 1])$
or ${\rm lk}(7) = C_7(6, 0, 1, x, [y, z, 9])$. But in both these cases ${\rm lk}(7)$ can not be completed. When ${\rm lk}(6) = C_7(0, 5, 9, 2, [8, d, 7])$, we see that $d = 1$ or 10. If $d = 1$ then $C_4(0, 1, 8, 6)\subseteq{\rm lk}(7)$. If $d = 10$ then ${\rm lk}(6) = C_7(0, 5, 9, 2, [8, 10,7 ])$ this implies ${\rm lk}(5) = C_7([11, 1, 9], 6, 0, 4, e)$ where $e\in\{7, 8, 10\}$. If $e = 7$ or $8$ then ${\rm lk}(7)$ has more than seven vertices. If $e = 10$ then ${\rm lk}(5) = C_7([9, 1, 11], 10, 4, 0, 6)$ this implies $C_5(6, 2, 1, 11, 5)\subseteq{\rm lk}(9)$. This is not possible.

\noindent{\bf Subcase 2.1:} When ${\rm lk}(6) = C_7(0, 5, 8, 2, [9, d, 7])$, we see again that $d = 1$ or 10. In case $d = 1$ we get $C_4(1, 0, 6, 9)\subseteq{\rm lk}(7)$. If $d = 10$, then ${\rm lk}(6) = C_7(0, 5, 8, 2, [9, 10, 7])$ and we get ${\rm lk}(5) = C_7([1, 11, 8], 6, 0, 4, e)$ or ${\rm lk}(5) = C_7([11, 1, 8], 6, 0, 4, e)$ for some $e \in V$. When ${\rm lk}(5) = C_7([11, 1, 8], 6, 0, 4, e)$, we see that $e \in \{7, 9, 10\}$. If $e  = 7$, then $C_(1, 0, 6, 9)\subseteq{\rm lk}(9)$. If $e = 9$, then ${\rm lk}(9)$ has more than seven vertices. So, $e = 10$. Then ${\rm lk}(8) = C_7([1, 11, 5], 6, 2, 3$, $y)$, where $y \in \{7, 9, 10\}$.

When $y = 7$, then ${\rm lk}(8) = C_7([1, 11, 5], 6, 2, 3, 7)$. Completing successively, we get  ${\rm lk}(5) = C_7([8, 1, 11], 10, 4, 0, 6)$, ${\rm lk}(1) = C_7([8, 5, 11], 9, 2, 0, 7)$, ${\rm lk}(7) = C_7([6, 9, 10], 3, 8, 1, 0)$, ${\rm lk}(3) = ([4, 0$, $2], 8, 7, 10, 11)$, ${\rm lk}(10) = C_7([7, 6, 9], 4, 5, 11, 3)$, ${\rm lk}(4) = C_7([3, 2, 0], 5, 10, 9, 11)$, ${\rm lk}(9) = C_7([6, 7, 10]$, $4, 11, 1, 2)$ and  ${\rm lk}(11) = C_7([5, 8, 1], 9, 4, 3, 10)$. This is {\bf $K_3$\,} of Example \ref{e1} with vertex 10 replaced by $u$ and vertex 11 replaced by $v$.

When $y = 9$, then ${\rm lk}(8) = C_7([1, 11, 5], 6, 2, 3, 9)$. Now, completing successively, we get ${\rm lk}(9) = C_7([6, 7, 10], 3, 8, 1, 2)$, ${\rm lk}(10) = C_7([9, 6, 7], 4, 5, 11, 3)$, ${\rm lk}(1) = C_7([8, 5, 11], 7, 0, 2, 9)$, ${\rm lk}(3) =C_7([2, 0$, $4], 11, 10, 9, 8)$, ${\rm lk}(11) = C_7([1, 8, 5], 10, 3, 4, 7)$, ${\rm lk}(7) = C_7([10, 9, 6], 0, 1, 11, 4)$ and ${\rm lk}(4) = C_7([3, 2, 0]$, $5, 10, 7, 11)$. This is {\bf $K_2$\,} of Example \ref{e1} with vertex 10 replaced by $u$ and vertex 11 replaced by $v$ for convenience.

\smallskip

In following lines the vertex set is retained as $\{0, 1, 2, \ldots, 10, 11\}$ for the sake of computational convenience. When we give isomorphism to $K_1$, $K_2$ or $K_3$, we mean isomorphism is given to the $K_1$ in Case 1, to $K_2$ in Subcase 2.1.2 and to $K_3$ in Subcase 2.1.1 above.

\noindent {\bf Subcase 2.2:} If ${\rm lk}(5) = C_7([1, 11, 8], 6, 0, 4, e)$ it is easy to see that $e\in\{7, 9, 10\}$. If $e = 10$, then ${\rm lk}(1)$ has more than seven vertices.

If $e = 7$, then completing successively we get ${\rm lk}(1)$ = $C_7([5, 8, 11], 9, 2, 0, 7)$, ${\rm lk}(7)$ = $C_7([10, 9, 6]$, $0, 1, 5, 4)$, ${\rm lk}(4)$ = $C_7([3, 2, 0], 5, 7, 10, 11)$, ${\rm lk}(9)$ = $C_7([6, 7, 10], 3, 11, 1, 2)$, ${\rm lk}(5)$ = $C_7([1, 11, 8], 6, 0, 4$, $7)$, ${\rm lk}(8)$ = $C_7([5, 1, 11], 10, 3, 2, 6)$,  ${\rm lk}(3)$ = $C_7([2, 0, 4], 11, 9, 10, 8)$, ${\rm lk}(11)$ = $C_7([1, 5, 8], 10, 4, 3, 9)$ and ${\rm lk}(10)$ = $C_7([7, 6, 9], 3, 8, 11, 4)$.
It is isomorphic to {\bf $K_2$} by the map $(0, 2)(3, 4)(5, 8)(7, 9)$.

\smallskip

If $e = 9$, then ${\rm lk}(5)$ = $C_7([8, 11, 1], 9, 4, 0, 6)$. Completing successively we get ${\rm lk}(1)$ = $C_7([5, 8, 11]$, $7, 0, 2, 9)$, ${\rm lk}(7)$ = $C_7([10, 9, 6], 0, 1, 11, 3)$, ${\rm lk}(3)$ = $C_7([2, 0, 4], 11, 7, 10, 8)$, ${\rm lk}(8)$ = $C_7([5, 1, 11], 10, 3$, $2, 6)$, ${\rm lk}(10)$ = $C_7([9, 6, 7], 3, 8, 11, 4)$, ${\rm lk}(4)$ = $C_7([3, 2, 0], 5, 9, 10, 11)$, ${\rm lk}(11)$ = $C_7([8, 5, 1], 7, 3, 4, 10)$ and  ${\rm lk}(9)$ = $C_7([6, 7, 10], 4, 5, 1, 2)$. It is isomorphic to {\bf $K_3$} by the map $(0, 6, 8, 3, 10, 11, 4, 9, 1)(2, 7, 5)$.

\smallskip

\noindent{\bf Case 3:} When $(a, b, c) = (8, 7, 6)$, we have ${\rm lk}(2) = C_7([3, 4, 0], 1, 6, 7, 8)$. This implies ${\rm lk}(6) = C_7(1, 2, 7, 0, [5, a, b])$ or ${\rm lk}(6) = C_7(5, 0, 7, 2, [1, a, b])$, for some $a, b \in V$.

If ${\rm lk}(6) = C_7(5, 0, 7, 2, [1, a, b])$, it is easy to see that  ${\rm lk}(6) = C_7([8, 9, 1], 2, 7, 0, 5)$, ${\rm lk}(6) = C_7([9, 8, 1], 2, 7, 0, 5)$ or ${\rm lk}(6) = C_7([9, 10, 1], 2, 7, 0, 5)$. When ${\rm lk}(6) = C_7([8, 9, 1], 2, 7, 0, 5)$ or ${\rm lk}(6) = C_7([9, 8, 1], 2, 7, 0, 5)$ then considering ${\rm lk}(7)$, it is easy to see that vertex $8$ or $1$ is in two quadrangles. This can not happen.  If ${\rm lk}(6) = C_7([9, 10, 1], 2, 7, 0, 5)$ then ${\rm lk}(7) = C_7([8, 5, 11], 1, 0, 6$, $2)$ or ${\rm lk}(7) =([8, 11, 5], 1, 0, 6, 2)$. In the first case we get $C_6(2, 0, 7, 10, 9, 6)\subseteq{\rm lk}(1)$. In the second case we see that ${\rm lk}(5)$ has more than seven vertices.

In the second case, we get ${\rm lk}(6) = C_7(1, 2, 7, 0, [5, 8, 9])$, ${\rm lk}(6) = C_7(1, 2, 7, 0, [5, 9, 8])$ or ${\rm lk}(6) = C_7(1, 2, 7, 0, [5, 9, 10])$. If ${\rm lk}(6) = C_7(1, 2, 7, 0, [5, 9, 8])$ then ${\rm lk}(7) = C_7([1, 10, 11], 8, 2, 6, 0)$. But then ${\rm lk}(8)$ has more than seven vertices. This is not possible. If ${\rm lk}(6) = C_7(1, 2, 7, 0, [5, 9, 10])$, then ${\rm lk}(1) = C_7([7, 8, 11], 10, 6, 2, 0)$. In this case we get $C_6 (6, 2, 8, 11, 1, 0)\subseteq{\rm lk}(7)$. This is not possible. When ${\rm lk}(6) = C_7(1, 2, 7, 0, [5, 8, 9])$ we get ${\rm lk}(7) = C_7(0, 6, 2, 8, [11, 10, 1])$ or ${\rm lk}(7) = C_7(0, 6, 2, 8, [10, 11, 1])$.

If ${\rm lk}(7) = C_7(0, 6, 2, 8, [11, 10, 1])$ then ${\rm lk}(1) = C_7([7, 11, 10], 9, 6, 2$, $0)$. This implies ${\rm lk}(8) = C_7([5, 6, 9], 11, 7, 2, 3)$ or ${\rm lk}(8) = C_7([5, 6, 9], 3, 2, 7, 11)$. If ${\rm lk}(8) = C_7([5, 6, 9], 3, 2, 7, 11)$ we have ${\rm lk}(9) = C_7(3, [8, 5, 6], 1, 10, b)$ where $b\in\{4, 11\}$. If $b = 4$ then ${\rm lk}(4) = C_7([0, 2, 3], 9, 10, 11, 5)$. This implies $C_6(5, 4, 10, 1, 7, 8)\subseteq {\rm lk}(11)$. Which is not allowed. If $b = 11$ then ${\rm lk}(11) = C_7([7, 1, 10], 9, 3, 5, 8)$, ${\rm lk}(3) = C_7([2, 0, 4], 5, 11, 9, 8)$. Then $C_4(5, 3, 2, 0) \subseteq {\rm lk}(4)$. This is not possible. So, ${\rm lk}(8) = C_7([5, 6, 9], 11, 7, 2, 3)$, completing successively we get ${\rm lk}(3) = C_7([4, 0, 2], 8, 5, 10$, $11)$, ${\rm lk}(5) = C_7([6, 9, 8], 3, 10, 4, 0)$, ${\rm lk}(9) = C_7([8, 5, 6], 1, 10, 4, 11)$, ${\rm lk}(4) = C_7([3, 2, 0], 5, 10, 9, 11)$, ${\rm lk}(10) = C_7([1, 7, 11], 3, 5, 4, 9)$, and ${\rm lk}(11) = C_7([7, 1, 10], 3, 4, 9, 8)$. This object is isomorphic to $K_1$ by the map $(0, 3)(1, 11)(2, 4)(5, 6, 9, 8)(7, 10)$.

If ${\rm lk}(7) = C_7(0, 6, 2, 8, [10, 11, 1])$ we have ${\rm lk}(6) = C_7([5, 8, 9], 1, 2, 7, 0)$, ${\rm lk}(8) = C_7([5, 6, 9], 10, 7$, $2, 3)$. This implies ${\rm lk}(10) = C_7([11, 1, 7], 8, 9, 3, 4)$ or ${\rm lk}(10) =  C_7([11, 1, 7], 8, 9, 4$, $3)$. If ${\rm lk}(10) = C_7([11, 1, 7], 8, 9, 3, 4)$ then ${\rm lk}(3) = C_7([4, 0, 2], 8, 5, 9, 10)$. But then $59$ form an edge which is not allowed since they are non edges in the quadrangle $[5, 6, 9, 8]$. If ${\rm lk}(10) = C_7([11, 1, 7], 8, 9, 4, 3)$, then we have ${\rm lk}(4) = C_7([0, 2, 3], 10, 9, 11, 5)$. Completing successively we get ${\rm lk}(5) = C_7([6, 9, 8], 3, 11$, $4, 0)$, ${\rm lk}(11) = C_7([1, 7, 10], 3, 5, 4, 9)$, ${\rm lk}(9) = C_7([6, 5, 8], 10, 4, 11, 1)$. This is isomorphic to $K_1$ by the map $(0, 3)(1, 11, 7, 10)(2, 4)(5, 6, 9, 8)$.

\smallskip

\noindent{\bf Case 4:} If $(a, b, c) = (8, 9, 6)$ then ${\rm lk}(2) = C_7([3, 4, 0], 1, 6, 9, 8)$. This implies ${\rm lk}(6) = C_7(0, 7, 9, 2, [1$, $a, 5])$  or ${\rm lk}(6) = C_7(0, 5, 1, 2, [9, a, 7])$ for some $a \in V$. If ${\rm lk}(6) = C_7(0, 7, 9, 2, [1, a, 5])$, it is easy to see that $a \in \{8, 10\}$.

If $a = 8$ then ${\rm lk}(6) = C_7(0, 7, 9, 2, [1, 8, 5])$ then  ${\rm lk}(9) = C_7([11, 10, 7], 6, 2, 8, b)$ where $b \in \{3, 4, 5\}$. If $b = 3$ then, considering ${\rm lk}(2)$ we see that  $C_3(9, 2, 3) \subseteq {\rm lk}(8)$. When $b = 4$ we have ${\rm lk}(8) = C_7([5, 6, 1], 3, 2, 9, 4)$ and ${\rm lk}(4) = C_7([3, 2, 0], 5, 8, 9, 11)$. In this case we get $C_(4, 0, 6, 1, 8)\subseteq{\rm lk}(5)$. When $b = 5$ we have ${\rm lk}(9) = C_7([7, 10, 11], 5, 8, 2, 6)$, ${\rm lk}(5) = C_7([6, 1, 8], 9, 11, 4, 0)$ and ${\rm lk}(8) = C_7([1, 6, 5], 9, 2, 3, c)$, for some  $c\in\{10, 11\}$. If $c = 10$ then $C_6(0, 1, 10, 11, 9, 6)\subseteq{\rm lk}(7)$. If $c = 11$, then ${\rm lk}(1) = C_7([6, 5, 8], 11, 7, 0, 2)$. Now,  considering${\rm lk}(9)$ we see that $7\,11$ are non edges because they for diagonal vertices of a quadrangle containing $9$ and in ${\rm lk}(1)$ they form an edge. A contradiction.
When $a = 10$ it is easy to see that ${\rm lk}(9) = C_7(6, 2, 8, a, [b, c, 7])$ or ${\rm lk}(9) = C_7(2, 6, 7, a, [b, c, 8])$. In both these cases ${\rm lk}(9)$ can not be completed. So, ${\rm lk}(6) \neq C_7(0, 7, 9, 2, [1, a, 5])$.

When ${\rm lk}(6) = C_7(0, 5, 1, 2, [9, a, 7])$ we have $a\in\{8, 10\}$. When $a = 8$ then ${\rm lk}(6) = C_7(0, 5, 1, 2$, $[9, 8, 7])$ ${\rm lk}(5) = C_7([11, 10, 1], 6, 0, 4, c)$ where  $c\in\{7, 8, 9\}$. In each of the respective cases, the links ${\rm lk}(7)$,  ${\rm lk}(8)$  and ${\rm lk}(9)$ have more than seven vertices respectively. When $a = 10$ then ${\rm lk}(5) = C_7([11, 8, 1], 6, 0, 4, b)$ or ${\rm lk}(5) = C_7([8, 11, 1], 6, 0, 4, b)$, for some $b \in V$.
\smallskip

{\bf Case 4.1 :} When ${\rm lk}(5) = C_7([11, 8, 1], 6, 0, 4, b)$ we have $b\in\{7, 9, 10\}$. In case  $b = 7$ or $b = 9$, the links of $7$ or $9$ have more than seven vertices. If $b = 10$ then ${\rm lk}(5) =  C_7([11, 8, 1], 6, 0, 4, 10)$. So, ${\rm lk}(1) = C_7([5, 11, 8], 7, 0, 2, 6)$ and ${\rm lk}(7) = C_7([10, 9, 6], 0, 1, 8, c)$ where $c\in\{3, 4, 11\}$. If $c = 4$ then ${\rm lk}(8)$ has more than seven vertices. If $c = 11$ then $C_4(1, 7, 11, 5)\subseteq{\rm lk}(8)$. So, $c = 3$. Then ${\rm lk}(7) = C_7([10, 9, 6], 0, 1, 8, 3)$. Now completing successively we get ${\rm lk}(8) = C_7([1, 5, 11], 9$, $2, 3, 7)$, ${\rm lk}(9) = C_7([10, 7, 6], 2, 8, 11, 4)$, ${\rm lk}(10) = C_7([7, 6, 9], 4, 5, 11, 3)$, ${\rm lk}(3) = C_7([4, 0, 2], 8, 7, 10$, $11)$, ${\rm lk}(4) = C_7([0, 2, 3], 11, 9, 10, 5)$ and ${\rm lk}(11) =C_7 ([5, 1, 8], 9, 4, 3, 10)$. This object is isomorphic to $K_2$ by the map $(0, 10, 8)$ $(1, 4, 9)$ $(2, 7, 5, 3, 6, 11)$.

\smallskip

{\bf Case 4.2 :} When ${\rm lk}(5) = C_7([8, 11, 1], 6, 0, 4, a)$ we see that $a\in\{9, 10\}$. If $a = 9$ we get ${\rm lk}(5) = C_7([1, 11, 8], 9, 4, 0, 6)$, ${\rm lk}(9) = C_7([6, 7, 10], 4, 5, 8, 2)$ and ${\rm lk}(4) = C_7([0, 2, 3], 11, 10, 9, 5)$. This implies ${\rm lk}(11) = C_7([1, 5, 8], 3, 4, 10, 7)$ or ${\rm lk}(11) = C_7([1, 5, 8], 10, 4, 3, 7)$. In the first case we get $C_(7, 11, 4, 9, 6)\subseteq {\rm lk}(6)$ which is not possible. So, ${\rm lk}(11) = C_7([1, 5, 8], 10, 4, 3, 7)$. Now, completing successively, we get ${\rm lk}(10) = C_7([7, 6, 9], 4, 11, 8, 3)$, ${\rm lk}(8) = C_7([5, 1, 11], 10, 3$, $2, 9)$, ${\rm lk}(3) = C_7([2, 0, 4], 11, 7, 10, 8)$ and ${\rm lk}(7) = C_7([6, 9, 10], 3, 11, 1, 0)$. This is isomorphic to $K_2$ by the map $(0, 11, 9)(1, 10, 2, 5, 7, 3, 8, 6, 4)$.

If $a = 10$ then ${\rm lk}(5) = C_7([1, 11, 8], 10, 4, 0, 6)$, ${\rm lk}(1) = C_7([5, 8, 11], 7, 0, 2, 6)$ this implies ${\rm lk}(8) = C_7([11, 1, 5], 10, 3, 2, 9)$ or ${\rm lk}(8) = C_7([11, 1, 5], 10, 9, 2, 3)$. If ${\rm lk}(8) = C_7([11, 1, 5],1 0, 9, 2, 3)$ then $C_5(2, 8, 10, 7, 6) \subseteq {\rm lk}(9)$. This is not possible. So, ${\rm lk}(8) = C_7([11, 1, 5], 10, 3, 2, 9)$ and we get ${\rm lk}(3) = C_7([4, 0, 2], 8, 10, b, a)$, for some $a, b \in V$. It is easy to see that $(a, b) \in \{(7, 5)$, $(7, 6)$, $(7, 11)$, $(9, 5)$, $(9, 11)$, $(11, 5)$, $(11, 7)$, $(11, 9)\}$. For the values $(a, b) \in \{(7, 5)$, $(7, 6)$, $(7, 11)$, $(9, 5)$, $(9, 11)$, $(11, 5)\}$ the link of one of the vertices $5, 6, 7$ or $9$ has more than seven vertices. Thus, $(a, b) \in \{(11, 7), (11, 9)\}$. Then ${\rm lk}(3) = C_7([4, 0, 2], 8, 10, 7, 11)$. Completing successively, we get  ${\rm lk}(7) = C_7([6, 9, 10], 3, 11$, $1, 0)$, ${\rm lk}(11) = C_7([8, 5, 1], 7, 3, 4, 9)$, ${\rm lk}(4) = C_7([0, 2, 3], 11, 9, 10, 5)$, ${\rm lk}(9) = C_7([6, 7, 10], 4, 11, 8, 2)$ and ${\rm lk}(10) = C_7([7, 6, 9], 4, 5, 8, 3)$.  This is isomorphic to $K_3$ by the map $(0, 11, 7, 3, 8, 6, 4, 1, 10, 2, 5, 9)$.

When $(a, b) = (11, 9)$ we have ${\rm lk}(3) = C_7([4, 0, 2], 8, 10, 9, 11)$. Now completing successively we get ${\rm lk}(9) =C_7([6, 7, 10], 3, 11, 8, 2)$, ${\rm lk}(11) = C_7([8, 5, 1], 7, 4, 3, 9)$, ${\rm lk}(7) = C_7([6, 9, 10], 4, 11, 1, 0)$, ${\rm lk}(4) = C_7([3, 2, 0],5, 10, 7, 11)$, ${\rm lk}(10) = C_7([7, 6, 9], 3, 8, 5, 4)$.  This is isomorphic to $K_2$ by the map $(0, 2, 3, 4)(1, 8, 11, 5)(6, 9, 10, 7)$.

\smallskip

\noindent{\bf Case 5:} When $(a, b, c) = (8, 9, 10)$ then ${\rm lk}(2) = C_7(1, [0, 4, 3], 8, 9, 10)$. This implies ${\rm lk}(1) = C_7(2, 0, 7, x, [y, z, 10])$ or ${\rm lk}(1) = C_7(10, 2, 0, 7, [x, y, z])$, for some $x, y, z \in V$.

\smallskip




\smallskip






When $(a, b, c) = (8, 9, 10)$ then ${\rm lk}(2) = C_7(1, [0, 4, 3], 8, 9, 10)$. This implies ${\rm lk}(1) = C_7(2, 0, 7$, $x, [y, z, 10])$ or ${\rm lk}(1) = C_7(10, 2, 0, 7, [x, y, z])$, for some $x, y, z \in V$.
\smallskip

\noindent {\bf Subcase 5.1:} When ${\rm lk}(1) = C_7(2, 0, 7, x, [y, z, 10])$, it is easy to see that $(x, y, z) \in \{(5, 6, 8)$, $(5, 6, 11)$, $(5, 8, 6)$, $(5, 8, 11)$, $(5, 11, 6)$, ($5, 11, 8)$, ($8, 5, 6)$, $(8, 5, 11)$, $(8, 6, 5)$, $(8, 6, 11)$, $(8, 11, 5)$, $(8, 11, 6)$, $(9, 5, 6)$, $(9, 5, 8)$, $(9, 5, 11)$, $(9, 6, 5)$, $(9, 6, 8)$, $(9, 6, 11)$, $(9, 8, 5)$, $(9, 8, 6)$, $(9, 8, 11)$, $(11, 5, 6)$, $(11, 5, 8)$, $(11, 6, 5)$, $(11, 6, 8)$, $(11, 8, 5)$, $(11, 8, 6)\}$.

\smallskip

\begin{claim}\label{cl3} There does not exist a SEM for the values $(x, y, z) \in \{(5, 6, 11)$, $(5, 8, 6)$, $(8, 6, 11)$, $(8, 11, 5)$, $(9, 5, 6)$, $(9, 8, 5)$, $(9, 8, 6)$, $(9, 8, 11)$, $(11, 5, 8)$, $(11, 6, 5)$, $(11, 6, 8)$, $(11, 8, 5)$, $(11, 8, 6)\}$.
\end{claim}

{\sc Proof of the Claim:}\ref{cl3} When $(x,y,z) = (5, 6, 11)$ we have ${\rm lk}(1) = C_7(2, 0, 7, 5, [10, 11, 6])$. This implies ${\rm lk}(6) = C_7([11, 10, 1], 5, 0, 7, a)$ where $a \in \{3, 4, 8, 9\}$. If $a = 3$ then ${\rm lk}(6) = C_7([1, 10, 11], 3, 7, 0, 5)$. This implies ${\rm lk}(3) = C_7([2, 0, 4], 7, 6, 11, 8)$ or ${\rm lk}(3) = C_7([2, 0, 4], 11, 6, 7, 8)$. In the first case, $i.e.$, ${\rm lk}(3) = C_7([2, 0, 4], 11, 6, 7, 8)$ then ${\rm lk}(7) = C_7([8, 9, 5], 1, 0, 6, 3)$. This implies $C_5(3, 2, 9, 5, 7)\subseteq{\rm lk}(8)$, which is not possible. If ${\rm lk}(3) = C_7([2, 0, 4], 7, 6, 11, 8)$ then considering ${\rm lk}(7)$, we see that $4$ lies in two different quadrangles. This is not possible. If $a = 4$ we get ${\rm lk}(6) = C_7([1, 10, 11], 4, 7, 0, 5)$. This implies ${\rm lk}(4) = C_7([0, 2, 3], 11, 6, 7, 5)$ or ${\rm lk}(4) = C_7([0, 2, 3], 7, 6, 11, 5)$. When ${\rm lk}(4) = C_7([0, 2, 3], 11, 6, 7, 5)$ then $C_5(5, 4, 6, 0, 1)\in{\rm lk}(7)$. If ${\rm lk}(4) = C_7([0, 2, 3], 7, 6, 11, 5)$ then considering ${\rm lk}(7)$, we get 3 lie in two different quadrangles. When $a = 8$ we have ${\rm lk}(6) = C_7([1, 10, 11], 8, 7, 0, 5)$. This implies ${\rm lk}(5) = C_7([7, 8, 9], 4, 0, 6, 1)$. Then $C_5(0, 6, 8, 9, 5)\subseteq{\rm lk}(7)$. When $a = 9$ we get ${\rm lk}(6) = C_7([1, 10, 11], 9, 7, 0, 5))$. This implies $C_6(0, 1, 5, 8, 9, 5)\subseteq{\rm lk}(7)$. 

\smallskip

When $(x, y, z) = (5, 8, 6)$ we have ${\rm lk}(1) = C_7(2, 0, 7, 5, [8, 6, 10])$. This implies ${\rm lk}(8)$ = $C_7(9, 2, 3$, $5,  [1, 10, 6])$ or $C_7(3, 2, 9, 5, [1, 10, 6])$. In both these cases there are more than seven vertices in ${\rm lk}(5)$. This is not possible. When $(x, y, z) = (5, 8, 11)$, we get ${\rm lk}(8) = C_7([1, 10, 11], 9, 2, 3, 5)$ or ${\rm lk}(8) = C_7([1, 10, 11], 3, 2, 9, 5)$. In this case ${\rm lk}(5)$ has more than seven vertices. When $(x, y, z) = (5, 11, 6)$ considering ${\rm lk}(5)$, we get $4$ or $6$ in two different quadrangles containing $5$. When $(x, y, z) = (5, 11, 8)$ we are in same situation as in previous case and hence this is also not possible. When $(x, y, z) = (8, 5, 6)$ then we get ${\rm lk}(8) = C_7([7, 11, 9], 2, 3, 5, 1)$ and ${\rm lk}(5) = C_7([1, 10, 6], 0, 4, 3, 8)$. This implies $C_4(8, 5, 4, 3)\in {\rm lk}(8)$. If $(x, y, z) = (8, 5, 11)$ then we have ${\rm lk}(8) = C_7([7, 6, 9], 2, 3, 5, 1)$. This implies $C_5(1, 0, 6, 9, 8)\subseteq{\rm lk}(7)$. 

\smallskip
If $(x, y, z) = (8, 6, 11)$ then ${\rm lk}(1) = C_7([10, 11, 6], 8, 7, 0, 2)$ and ${\rm lk}(8) = C_7([7, 5, 9], 2, 3, 6, 1)$. This implies ${\rm lk}(9) = C_7([5, 7, 8], 2, 10, a, b)$ where possible values of $(a, b)$ are $\{(3, 4)$, $(4, 3)$, $(4, 11)$, $(11, 4)\}$. If $(a, b) = (3, 4)$ then ${\rm lk}(9) = C_7([5, 7, 8], 2, 10, 3, 4)$. Then ${\rm lk}(4)$ can not be completed. If $(a, b) = (4, 3)$ then ${\rm lk}(9) = C_7([5, 7, 8], 2, 10, 4, 3)$. This implies ${\rm lk}(4) = C_7(10, 9, [3, 2, 0], 5, c)$ where $c \in \{6, 7, 8, 11\}$. For the successive values $c = 6, 7, 8, 11$ we see that ${\rm lk}(6) , {\rm lk}(7), {\rm lk}(8)$ and ${\rm lk}(5)$ have more than seven vertices. If $(a, b) = (4, 11)$ we have ${\rm lk}(4) = C_7([0, 2, 3], 11, 9, 10, 5)$ or ${\rm lk}(4) = C_7([0, 2, 3], 10, 9, 11, 5)$. When ${\rm lk}(4) = C_7([0, 2, 3], 11, 9, 10, 5)$ we have ${\rm lk}(10) = C_7([1, 6, 11], 5, 4, 9$, $2)$. In this case ${\rm lk}(5)$ has more than seven vertices. When ${\rm lk}(4) = C_7([0, 2, 3], 10, 9, 11, 5)$, we get $C_3(4, 9, 5) \subseteq {\rm lk}(11)$. If $(a, b) = (11, 4)$ we have ${\rm lk}(9) = C_7([8, 7, 5], 4, 11, 10, 2)$. Then, we get $C_5(2, 9, 11, 6, 1) \in {\rm lk}(10)$. So, $(x, y, z) \neq (8, 6, 11)$.

\smallskip

When $(x, y, z) = (8, 11, 5)$ we get ${\rm lk}(7) = C_7(0, 1, 8, p, [q, r, 6])$ or ${\rm lk}(7) = C_7(1, 0, 6, p, [q, r, 8])$ for some $p, q, r \in V(K)$. In both these cases no suitable value of $q$ and $r$ can be found in $V(K)$. Hence $(x, y, z) \neq (8, 11, 5)$. When $(x, y, z) = (9, 5, 6)$, we have ${\rm lk}(9) = C_7([7, 11, 8], 2, 10, 5, 1)$. This implies ${\rm lk}(5) = C_7(9, [1, 10, 6], 0, 4, a)$, where $a \in \{3, 11\}$. If $a = 3$ then $C_3(3, 0, 5) \subseteq {\rm lk}(4)$. Also if $a = 11$ then $9\, 11$ is and edge. But, considering ${\rm lk}(9)$, we see that this is a non - edge. When $(x, y, z) = (9, 5, 8)$ or $(9, 5, 11)$ then either 8 or 10 are in two different quadrangles. Which is not possible. When $(x, y, z) = (9, 6, 5)$ then ${\rm lk}(9) = C_7([7, 11, 8], 2, 10, 6, 1)$. Considering ${\rm lk}(1)$ we see that $6\,10$ is non-edge. When $(x, y, z) = (9, 6, 8)$ then for the same reason as in previous case we see that this is also not possible. When $(x, y, z) = (9, 6, 11)$ then we have ${\rm lk}(9) = C_7([7, 5, 8], 2, 10, 6, 1)$ but again $6\,10$ is a non-edge in the ${\rm lk}(1)$. So, $(x, y, z) \neq (8, 11, 5)$ or $(9, 5, 6)$.

\smallskip

When $(x, y, z) = (9, 8, 5)$ then ${\rm lk}(1) = C_7([10, 5, 8], 9, 7, 0, 2)$ and ${\rm lk}(9) = C_7([7, 6, 11], 10, 2, 8, 1)$. This implies ${\rm lk}(10) = C_7([5, 8, 1], 2, 9, 11, c)$ where $c \in \{3, 4, 6, 7\}$. If $c = 3$ then ${\rm lk}(10) = C_7([1, 8, 5], 3, 11, 9, 2)$ . This implies ${\rm lk}(5) = C_7([8, 1, 10], 3, 6, 0, 4)$, ${\rm lk}(8) = C_7([5, 10, 1], 9, 2, 3, 4)$ then we get $C_3(3, 8, 5, 0) \subseteq {\rm lk}(4)$. If $c = 4$ then ${\rm lk}(10) = C_7([1, 8, 5], 4, 11, 9, 2)$ this implies ${\rm lk}(8) = C_7([5, 10, 1], 9, 2, 3, d)$ where$d \in \{4, 6\}$. If $d = 6$ then ${\rm lk}(8) = C_7([5, 10, 1], 9, 2, 3, 6)$ this implies ${\rm lk}(6) = C_7([11, 9, 7], 0, 5, 8, 3)$ then $C_5(0, 1, 9, 11, 6)\subseteq{\rm lk}(7)$. If $d = 4$ then ${\rm lk}(8) = C_7([5, 10, 1], 9, 2, 3, 4)$ then we get $C_3(8, 3, 0, 5) \subseteq {\rm lk}(4)$. If $c = 6$ we get ${\rm lk}(10) = C_7([1, 8, 5], 6, 11, 9$, $2)$. Then $C_6(0, 5, 10, 11, 9, 7)\subseteq {\rm lk}(6)$. If $c = 7$ then ${\rm lk}(10) = C_7([1, 8, 5], 7, 11, 9, 2)$. We get seven vertices in ${\rm lk}(7)$. This is not possible. So, $(x, y, z) \neq (9, 8, 5)$.

\smallskip

When $(x, y, z) = (9, 8, 6)$ then ${\rm lk}(1) = C_7([10, 6, 8], 9, 7, 0, 2)$ this implies ${\rm lk}(9) = C_7([5, 11, 7], 1$, $8, 2, 10)$ or ${\rm lk}(9) = C_7([11, 5, 7], 1, 8, 2, 10)$. If ${\rm lk}(9) = C_7([5, 11, 7], 1, 8, 2, 10)$ then we get ${\rm lk}(5) = C_7([9, 7, 11], 6, 0, 4, 10)$ or ${\rm lk}(5) = C_7([9, 7, 11], 4, 0, 6, 10)$. When ${\rm lk}(5) = C_7([9, 7, 11], 6, 0, 4, 10)$ we have ${\rm lk}(10) = C_7([1, 8, 6], 4, 5, 9, 2)$. In this case we find more than seven vertices in ${\rm lk}(6)$. This is not possible. If ${\rm lk}(5) = C_7([9, 7, 11], 4, 0, 6, 10)$ then $C_6(2, 9, 5, 6, 8, 1)\subseteq{\rm lk}(10)$. If ${\rm lk}(9) = C_7([11, 5, 7], 1, 8, 2, 10)$ we have ${\rm lk}(7) = C_7([5, 11, 9], 1, 0, 6, b)$, for some $b \in \{3, 4, 8, 10\}$. If $b = 3$,then ${\rm lk}(3) = C_7([2, 0, 4], 6, 7, 5, 8)$. This implies there are more than seven vertices in ${\rm lk}(5)$. This is not possible. If $b = 4$ then ${\rm lk}(7) = C_7([5, 11, 9], 1, 0, 6, 4)$. So, ${\rm lk}(4) = C_7([3, 2, 0], 5, 7, 6, c)$ for some $c \in \{8, 10, 11\}$. When $c = 8$, we have ${\rm lk}(4) = C_7([3, 2, 0], 5, 7, 6, 8)$ and ${\rm lk}(8) = C_7([1, 10, 6], 4, 3, 2, 9)$. So, ${\rm lk}(6) = C_7([8, 1, 10], 5, 0, 7, 4)$. Then ${\rm lk}(5) = C_7([7, 9, 11], 10, 6, 0, 4)$. Then $C_(10, 9, 7, 5)\subseteq{\rm lk}(11)$. This is not possible. When $c = 10$, we get seven vertices in ${\rm lk}(10)$. This is not possible. When $c = 11$, we get seven vertices in ${\rm lk}(6)$. If $b = 8$ or $10$ then ${\rm lk}(b)$ has more than seven vertices. This is not possible. Therefore $(x, y, z) \neq (9, 8, 6)$.

\smallskip

When $(x, y, z) = (9, 8, 11)$ we have ${\rm lk}(1) = C_7([10, 11, 8], 9, 7, 0, 2)$. This implies ${\rm lk}(9) = C_7([7, 6, 5], 10, 2, 8, 1)$. Then $C_4(7, 0, 5, 9)\subseteq{\rm lk}(6)$. When $(x, y, z) = (11, 5, 6)$, we have ${\rm lk}(1) = C_7([10, 11, 8], 9, 7, 0, 2)$. This implies ${\rm lk}(5) = C_7(4, 0, [6, 10, 1], 11, b)$ for some $b \in \{8, 9\}$. If $b = 8$ then ${\rm lk}(5) = C_7(4, 0, [6, 10, 1], 11, 8)$ this implies ${\rm lk}(8) = C_7([11, 7, 9], 2, 3, 4, 5)$ but $34$ and $38$ are adjacent edges in ${\rm lk}(2)$, while $348$ is a face in ${\rm lk}(8) = C_7([11, 7, 9], 2, 3, 4, 5)$. If $b = 9$ then ${\rm lk}(5) = C_7(4, 0, [6, 10, 1], 11, 9)$ and ${\rm lk}(9) = C_7([11, 7, 8], 2, 10, 4, 5)$. This implies ${\rm lk}(4) = C_7([3, 2, 0], 5, 9, 10, a)$ for some $a \in \{6, 7, 8, 11\}$. If $a = 6$ then ${\rm lk}(4) = C_7([3, 2, 0], 5, 9, 10, 6)$ and ${\rm lk}(6) = C_7([5, 1, 10], 4, 3, 7, 0)$. Then $C_6(2, 9, 4, 6, 5, 1)\subseteq {\rm lk}(10)$. If $a = 7$ then ${\rm lk}(4) = C_7([3, 2, 0], 5, 9, 10, 7)$. Then we get more than seven vertices in link of $7$. This is not possible. If $a = 8$ then ${\rm lk}(4) = C_7([3, 2, 0], 5, 9, 10, 8)$. This implies ${\rm lk}(8) = C_7([9, 11, 7], 10, 4, 3, 2)$. Then we get more than seven vertices in link of $10$. This is not possible. If $a = 11$ we have ${\rm lk}(4) = C_7([3, 2, 0], 5, 9, 10, 11)$. Then ${\rm lk}(11)$ has more than seven vertices. This is not possible. So, $(x, y, z) \neq (9, 8, 11)$.

\smallskip

When $(x, y, z) = (11, 5, 8)$ then ${\rm lk}(1) = C_7([10, 8, 5], 11, 7, 0, 2)$ this implies ${\rm lk}(8) = C_7([10, 1, 5]$, $3, 2, 9, b)$ for some $b \in \{4, 6, 7, 11\}$. If $b = 4$ then ${\rm lk}(8) = C_7([10, 1, 5], 3, 2, 9, 4)$. This implies ${\rm lk}(4) = C_7([0, 2, 3], 10, 8, 9, 5)$. Then ${\rm lk}(10) = C_7([8, 5, 1], 2, 9, 3, 4)$. In this case, considering ${\rm lk}(9)$ we see that $3$ or $5$ lie in two different quadrangles. If $b = 6$ then ${\rm lk}(8) = C_7([10, 1, 5], 3, 2, 9, 6)$. Considering ${\rm lk}(6)$ we see that $5$ lies in two different quadrangles. If $b = 7$ then ${\rm lk}(8) = C_7([10, 1, 5], 3, 2, 9, 6)$. This implies ${\rm lk}(7)$ has more than seven vertices. This is not possible. If $b = 11$ then ${\rm lk}(8) = C_7([10, 1, 5], 3, 2, 9, 6)$. Then $5$ or $10$ lie in two different quadrangles. Thus $(x, y, z) \neq (11, 5, 8)$.

\smallskip

If $(x, y, z) = (11, 6, 5)$ then ${\rm lk}(1) = C_7([6, 5, 10], 2, 0, 7, 11)$. This implies ${ \rm lk}(6) = C_7([11, [1, 10, 5]$, $0, 7, b)$ for some $b \in \{3, 4, 8, 9\}$. If $b = 3$ then ${\rm lk}(6) = C_7([11, [1, 10, 5], 0, 7, 3)$. This implies ${ \rm lk}(7) = C_7([11, 8, 9], 3, 6, 0, 1)$ or ${\rm lk}(7) = C_7([11, 9, 8], 3, 6, 0, 1)$. In the first case when ${ \rm lk}(7) = C_7([11, 8, 9], 3, 6, 0, 1)$ then ${\rm lk}(3)$ has more than seven vertices. In the second case ${\rm lk}(7) = C_7([11, 9, 8]$, $3, 6, 0, 1)$, $C_5(2, 9, 11, 7, 3)\subseteq{ \rm lk}(8)$. If $b = 4$ then ${\rm lk}(6) = C_7([11, [1, 10, 5], 0, 7, 4)$. This implies ${\rm lk}(7) = C_7([11, 9, 8], 4, 6, 0, 1)$ or ${\rm lk}(7) = C_7([11, 8, 9], 4, 6, 0, 1)$. In both these cases ${\rm lk}(4)$ has more than seven vertices. If $b = 8$ then ${\rm lk}(6) = C_7([11, [1, 10, 5], 0, 7, 8)$. So, ${\rm lk}(8) = C_7([11, k, 9], 2, 3, 7, 6)$ or ${\rm lk}(8) = C_7([7, k, 9], 2, 3, 11, 6)$. In both these cases the values of $k$ could not be found such that the links can be completed. If $b = 9$ then ${\rm lk}(6) = C_7([11, [1, 10, 5], 0, 7, 9)$. So, ${\rm lk}(9) = C_7([11, c, 8], 2, 10, 7, 6)$ or ${\rm lk}(9) = C_7([7, c, 8], 2, 10, 11, 6)$. As in the previous case, there is no value of $c$ such that the links can be completed. So, $(x, y, z) \neq (11, 6, 5)$.

\smallskip

When $(x, y, z) = (11, 6, 8)$ then ${\rm lk}(1) = C_7([6, 8, 10], 2, 0, 7, 11)$. This implies ${\rm lk}(6) = C_7([1, 10, 8]$, $7, 0, 5, 11)$ or ${\rm lk}(6) = C_7([1, 10, 8], 5, 0, 7, 11)$. If ${\rm lk}(6) = C_7([1, 10, 8], 7, 0, 5, 11)$ then ${\rm lk}(7) = C_7([11, 5, 9], 8, 6, 0, 1)$ and ${\rm lk}(9) = C_7([5, 11, 7], 8, 2, 10, 4)$. This implies $C_6(4, 0, 6, 11, 7, 9)\subseteq{\rm lk}(5)$. When ${\rm lk}(6) = C_7([1, 10, 8], 5, 0, 7, 11)$, we get $C_3(1, 0, 6, 11) \subseteq {\rm lk}(7)$. This is not possible. So, $(x, y, z) \neq (11, 6, 8)$.

\smallskip

When $(x, y, z) = (11, 8, 5)$ we have ${\rm lk}(1) = C_7([8, 5, 10], 2, 0, 7, 11)$. This implies ${\rm lk}(8) = C_7([1, 10, 5], 9, 2, 3, 11)$ or ${\rm lk}(8) = C_7([1, 10, 5], 3, 2, 9, 11)$. When ${\rm lk}(8) = C_7([1, 10, 5], 9, 2, 3, 11)$ we get ${\rm lk}(5) = C_7([8, 1, 10], 6, 0, 4, 9)$ or ${\rm lk}(5) = C_7([8, 1, 10], 4, 0, 6, 9)$. If ${\rm lk}(5) = C_7([8, 1, 10], 6, 0, 4, 9)$ then ${\rm lk}(9) = C_7(4, 5, 8, 2, [10, a, b])$. It is easy to see that $(a, b)\in\{(6, 7)$, $(6, 11)$, $(7, 6)$, $(7 ,11)$, $(11, 6)$, $(11, 7)\}$. When $(a, b) = (6, 7)$ then ${\rm lk}(9) = C_7([10, 6, 7], 4, 5, 8, 2)$. This implies ${\rm lk}(7) = C_7([9, 10, 6], 0, 1, 11, 4)$ and ${\rm lk}(4) = C_7([3, 2, 0], 5, 9, 7, 0)$. Then $C_4(3, 8, 1, 7, 4)\subseteq {\rm lk}(11)$, which is not possible. When $(a, b) = (6, 11)$ then ${\rm lk}(9) = C_7([10, 6, 11], 4, 5, 8, 2)$. This implies ${\rm lk}(4) = C_7([3, 2, 0], 5, 9, 11, c)$, where $c \in \{6, 7, 10\}$. If $c = 6$ or $10$ then link of $c$ has more than seven vertices in each case. If $c = 7$ then $3$ or $6$ appear in two quadrangles. When $(a, b) = (7, 6)$ we have ${\rm lk}(9) = C_7([10, 7, 6], 4, 5, 8, 2)$. In this case ${\rm lk}(4) = C_7([3, 2, 0], 5, 9, 6, d)$ for some $d \in \{10, 11\}$. In case $d = 10$ its link has more than seven vertices. If $d = 11$ then $C_5(8, 2, 0, 4, 11)\subseteq{\rm lk}(3)$. A contradiction. When $(a, b) = (7, 11)$ we have ${\rm lk}(9) = C_7([10, 7, 11], 4, 5, 8, 2)$. Then ${\rm lk}(9)$ is not possible because 3 or 4 will appear in two squares. When $(a, b) = (11, 6)$ we have ${\rm lk}(9) = C_7([10, 11, 6], 4, 5, 8, 2)$. This implies ${\rm lk}(4) = C_7([3, 2, 0], 5, 9, 11, d)$, for some $d$. But then $9$ and $11$ form an edge in ${\rm lk}(4)$. This is a contradiction. When $(a, b) = (11, 7)$ we have ${\rm lk}(9) = C_7([10, 11, 7], 4, 5, 8, 2)$. This implies ${\rm lk}(4) = C_7([3, 2, 0], 5, 9, 7, d)$ where $d \in\{10, 11\}$. If $d = 10$, ${\rm lk}(10)$ has more than seven vertices. If $d = 11$ then $C_5(8, 2, 0, 4, 11)\in {\rm lk}(3)$. Thus ${\rm lk}(5) \neq([8, 1, 10], 6, 0, 4, 9)$. Now, let ${\rm lk}(5) = C_7([8, 1, 10], 4, 0, 6, 9)$. This implies ${\rm lk}(9) = C_7(10, 2, 8, 5, [6, a, 11])$ for some $a$. But we see that no value of $a$ exist such that the SEM can be completed. Thus, ${\rm lk}(8) \neq([1, 10, 5], 9, 2, 3, 11)$. When ${\rm lk}(8) = C_7([1, 10, 5], 3, 2, 9, 11)$ we get ${\rm lk}(5) = C_7([8, 1, 10], 4, 0, 6, 3)$, and ${\rm lk}(10) = C_7(4, [5, 8, 1], 2, 9, b)$ where $b \in \{3, 7, 11\}$. If $b = 3$ or $7$ then link of $b$ has more than seven vertices. In case $b = 11$ then vertex $1$ appears in two quadrangles. So, $(x, y, z) \neq (11, 8, 5)$.

\smallskip

When $(x, y, z) = (11, 8, 6)$ we have ${\rm lk}(1) = C_7([8, 6, 10], 2, 0, 7, 11)$. This implies ${\rm lk}(6) = C_7([10, 1, 8], 5, 0, 7, b)$ or ${\rm lk}(6) = C_7([8, 1, 10], 7, 0, 5, b)$, for some $b \in V(K)$. If ${\rm lk}(6) = C_7([10, 1, 8], 5$, $0, 7, b)$, it is easy to see that $b \in \{3, 4, 9, 11\}$. If $b = 3$ then ${\rm lk}(7) =([11, 5, 9], 3, 6, 0, 1)$ or ${\rm lk}(7) = C_7([11, 9, 5], 3, 6, 0, 1)$. If ${\rm lk}(7) = ([11, 5, 9], 3, 6, 0, 1)$ then we get ${\rm lk}(9) = C_7([7, 11, 5], 8, 2, 10, 3)$. Then ${\rm lk}(8)$
has more than seven vertices. So, ${\rm lk}(7) = C_7([11, 9, 5], 3, 6, 0, 1)$. Then ${\rm lk}(3)$ has more than seven vertices. Which is not allowed. If $b = 4$ then ${\rm lk}(6) = C_7([8, 1, 10], 4, 7, 0, 5)$. This implies ${\rm lk}(4) = C_7([0, 2, 3], 10, 6, 7, 5)$ or ${\rm lk}(4) = C_7([0, 2, 3], 7, 6, 10, 5)$. In the first case, we have ${\rm lk}(10) = C_7([1, 8, 6], 4, 3, 9, 2)$ and ${\rm lk}(5) = C_7([7, 11, 9], 8, 6, 0, 4)$. Then ${\rm lk}(8)$ has more than seven vertices. In the second case, $3$ appears in two quadrangles. If $b = 9$ then ${\rm lk}(6) = C_7([8, 1, 10], 9, 7, 0, 5)$. This implies ${\rm lk}(9) = C_7([7, 11, 5], 8, 2, 10, 6)$, ${\rm lk}(5) = C_7([9, 7, 11], 4, 0, 6, 8)$ and ${\rm lk}(8) = C_7([6, 10, 1], 11, 2, 9, 5)$. Then ${\rm lk}(2)$ has more than 7 vertices. If $b = 11$ then ${\rm lk}(6) = C_7([8, 1, 10], 11, 7, 0, 5)$.
This implies $C_(11, 1, 0, 6) \in {\rm lk}(7)$. This is not possible. So, ${\rm lk}(6) \neq ([10, 1, 8], 5, 0, 7, b)$. If ${\rm lk}(6) = C_7([8, 1, 10], 7, 0, 5, b)$ it is easy to see that $b \in \{3, 4, 9, 11\}$. If $b = 3$ then ${\rm lk}(6) = C_7([8, 1, 10], 7, 0, 5, 3)$ and ${\rm lk}(8) = C_7([6, 10, 1], 11,
 9, 2, 3)$. This implies ${\rm lk}(10) = C_7(7, [6, 8, 1], 2, 9, c)$ for some $c \in \{3, 5, 11\}$. If $c = 3$ or $5$, then ${\rm lk}(c)$ has more than seven vertices. If $c = 11$ then $C_(1, 8, 9, 10, 7) \in {\rm lk}(11)$. This is not possible. If $b = 4$ then $C_(4, 5, 6)\in{\rm lk}(6)$. If $b = 9$. So, ${\rm lk}(6) = C_7([8, 1, 10], 7, 0, 5, 9)$ and ${\rm lk}(8) = C_7([6, 10, 1], 11, 3, 2, 9)$. This implies ${\rm lk}(9) = C_7([5, 7, 11], 10, 2, 8, 6)$ or ${\rm lk}(9) = C_7([5, 11, 7], 10, 2, 8, 6)$. If ${\rm lk}(9) = C_7([5, 7, 11], 10, 2, 8, 6)$ then we have ${\rm lk}(5) = C_7([7, 11, 9], 6, 0, 4, b)$ where $b\in\{8, 10\}$. In both these cases we get more than seven vertices in ${\rm lk}(c)$. If ${\rm lk}(9) = C_7([5, 11, 7], 10, 2, 8, 6)$ then ${\rm lk}(7) = C_7([9, 5, 11], 1, 0, 6, 10)$ and ${\rm lk}(5) = C_7([9, 7, 11], 10, 4, 0, 6)$. Then ${\rm lk}(10)$ has more than seven vertices. If $b = 11$ then we get ${\rm lk}(6) = C_7([8, 1, 10], 7, 0, 5, 11)$. Then $C_3(11, 1, 10, 6) \subseteq {\rm lk}(8)$. This is not possible. So, $(x, y, z) \neq (11, 8, 6)$. This Proves the Claim \ref{cl3}.

\smallskip
When $(x, y, z) = (5, 6, 8)$ we have ${\rm lk}(6) = C_7([8, 10, 1], 5, 0, 7, w)$ where $w \in \{3, 4, 9, 11\}$. If $w = 4$ we have ${\rm lk}(6) = C_7([8, 10, 1], 5, 0, 7, 4)$. This implies ${\rm lk}(5) = C_7(4, 0, 6, 1, [7, p, q])$. It is easy to see that $\{p, q\} = \{9, 11\}$. But in both the cases $(p, q) \in \{(9, 11), (11, 9)\}$, ${\rm lk}(4)$ has more than seven vertices. This is not possible. When $w = 9$ we have ${\rm lk}(6) = C_7([1, 10, 8],9, 7, 0, 5)$. This implies ${\rm lk}(9) = C_7(10, 2, 8, 6, [7, p, q])$, where $\{p, q\}\in \{5, 11\}$. When $(p, q) = (5, 11)$ we have ${\rm lk}(9) = C_7(10, 2, 8, 6, [7, 11, 5])$. Which is not possible, since $57$ is an edge in ${\rm lk(1)}$. When $(p, q) = (11, 5)$ we have ${\rm lk}(9) = C_7(10, 2, 8, 6, [7, 5, 11])$. Then $C_6(0, 1, 5, 11, 9, 6)\subseteq{\rm lk}(7)$. When $w = 11$ we have ${\rm lk}(6) = C_7(11, 7, 0, 5, [1, 10, 8])$. This implies ${\rm lk}(8) = C_7(3, 2, 9, 11, [6, 1, 10])$ or ${\rm lk}(8) = C_7(9, 2, 3, 11, [6, 1, 10])$. If ${\rm lk}(8) = C_7(9, 2, 3, 11, [6, 1, 10])$ then $C_3(8, 2, 10) \in {\rm lk}(9)$. If ${\rm lk}(8) = C_7(3, 2, 9, 11, [6, 1, 10])$ then ${\rm lk}(10) = C_7(3, [8, 6, 1], 2, 9, p)$ where $p \in \{4, 5, 7, 11\}$. When $p = 4$ then, $C_5(4, 10, 8, 2, 0)\subseteq{\rm lk}(3)$. If $p = 5$ then ${\rm lk}(5)$ has more than seven vertices. If $p = 7$ then ${\rm lk}(7)$ has more than seven vertices. If $p = 11$ then we see that either $3$ or $6$ lie in two quadrangles, whereas they are disjoint in $K$. So, $w = 3$. In this case we get ${\rm lk}(5) = C_7(4, 0, 6, 1, [7, p, q])$. It is easy to see that $\{p, q\} \in \{9, 11\}$. When $(p, q) = (9, 11)$ then ${\rm lk}(5) = C_7(4, 0, 6, 1, [7, 9, 11])$, ${\rm lk}(7) = C_7(3, 6, 0, 1, [5, 11, 9])$. This implies ${\rm lk}(9) = C_7([7, 5, 11], 8, 2, 10, 3)$ or ${\rm lk}(9) = C_7([7, 5, 11], 10, 2, 8, 3)$. When ${\rm lk}(9) = C_7([7, 5, 11], 8, 2, 10, 3)$, ${\rm lk}(3)$ has more than seven vertices. When ${\rm lk}(9) = C_7([7, 5, 11], 10, 2, 8, 3)$, we get $C_3(2, 3, 9)\subseteq {\rm lk}(8)$. When $(p, q) = (11, 9)$, ${\rm lk}(5) = C_7([7, 11, 9], 4, 0, 6, 1)$ and ${\rm lk}(7) = C_7([5, 9, 11], 3, 6, 0, 1)$. This implies ${\rm lk}(9) = C_7([5, 7, 11], 10, 2, 8, 4)$ or ${\rm lk}(9) = C_7([5, 7, 11], 2, 8, 10, 4)$. In the first case, $i.e.$, ${\rm lk}(9) = C_7([5, 7, 11], 10, 2, 8, 4)$ we have ${\rm lk}(8) = C_7([10, 1, 6], 3, 2, 9, 4)$, ${\rm lk}(10) = C_7([8, 6, 1], 2, 9, 11$, $4)$, ${\rm lk}(3) = C_7([2, 0, 4], 11, 7, 6, 8)$. Then $C_6(4, 3, 7, 5, 9, 10)\subseteq{\rm lk}(11)$. If ${\rm lk}(9) = C_7([5, 7, 11], 8, 2, 10$, $4)$ then completing successively we get ${\rm lk}(8) = C_7([10, 1, 6], 3, 2, 9, 11)$, ${\rm lk}(10) = C_7([8, 6, 1], 2, 9, 4$, $11)$, ${\rm lk}(11) = C_ 7([9, 5, 7], 3, 4, 10, 8)$, ${\rm lk}(3) = C_7([4, 0, 2], 8, 6, 7, 11)$, and ${\rm lk}(4) = C_7([3, 2, 0], 5, 9, 10$, $11)$. This is isomorphic to $K_1$ by the map $(0, 9, 2, 8, 1, 11)(3, 5)(4, 6, 10, 7)$.

\smallskip

When $(x, y, z) = (8, 6, 5)$ then ${\rm lk}(8) = C_7([7, 11, 9], 2, 3, 6, 1)$. This implies,${\rm lk}(6) = C_7([1, 10, 5]$, $0, 7, 3, 8)$, ${\rm lk}(7) = C_7([8, 9, 11], 3, 6, 0, 1)$, ${\rm lk}(3)= C_7([4, 0, 2], 8, 6, 7, 11)$, ${\rm lk}(10) =([5, 6, 1], 2, 9, a, b)$ for some $a, b \in \{4, 11\}$. If $(a, b) = (11, 4)$ we have ${\rm lk}(10) = C_7([1, 6, 5], 4, 11, 9, 2)$. This implies $C_6(3, 2, 0, 5, 10, 11)\subseteq{\rm lk}(4)$. So $(a, b) = (4, 11)$. This implies ${\rm lk}(10) = C_7([1, 6, 5], 11, 4, 9, 2)$. Now, completing successively we get ${\rm lk}(4) = C_7([3, 2, 0], 5, 9, 10, 11)$, ${\rm lk}(9) = C_7([8, 7, 11], 5, 4, 10, 2)$, ${\rm lk}(11) = C_7([9, 8, 7], 3, 4, 10, 5)$, ${\rm lk}(5) = C_7([10, 1, 6], 0, 4, 9, 11)$. It is isomorphic to $K_1$ by the map $(0, 8)(2, 5, 11, 3, 6, 7)(4, 9)$.

\smallskip

When $(x, y, z) = (8, 11, 6)$ we have ${\rm lk}(1) = C_7([10, 6, 11], 8, 7, 0, 2)$ and ${\rm lk}(8) = C_7([7, 5, 9], 2, 3$, $11, 1)$. This implies${\rm lk}(7) = C_7([5, 9, 8], 1, 0, 6, b)$ for some $b \in \{4, 10, 11\}$. When $b = 10$ we get ${\rm lk}(7) = C_7([8, 9, 5], 10, 6, 0, 1)$. We get $9$ or $1$ in two different quadrangles. If $b = 11$, then ${\rm lk}(7) = C_7([8, 9, 5], 11, 6, 0, 1)$ and ${\rm lk}(11) = C_7([1, 10, 6], 7, 5, 3, 8)$. Then we get more than seven vertices for in link of $5$. This is not possible. If $b = 4$, then ${\rm lk}(7) = C_7([8, 9, 5], 4, 6, 0, 1)$. This implies ${\rm lk}(5) = C_7([9, 8, 7], 4, 0, 6, c)$ where $c \in \{10, 11\}$. If $c = 10$ then ${\rm lk}(10) = C_7([1, 11, 6], 5, 9, 2)$. But then $C_6(9, 5, 6, 11, 1, 2)\subseteq{\rm lk}(10)$. Hence $c = 11$. So, ${\rm lk}(5) = C_7([7, 8, 9], 11, 6, 0, 4)$. Now completing successively we get ${\rm lk}(6) = C_7([10, 1, 11], 5, 0, 7, 4)$, ${\rm lk}(4) = C_7([0, 2, 3], 10, 6, 7, 5)$, ${\rm lk}(11) = C_7([6, 10, 1], 8, 3, 9, 5)$, ${\rm lk}(3) = C_7([2, 0, 4], 10, 9, 11, 8)$, ${\rm lk}(9)= C_7([5, 7, 8], 2, 10, 3, 11)$, ${\rm lk}(10) = C_7([1, 11, 6], 4, 3, 9, 2)$. This object is isomorphic to $K_1$ by the map $(0, 4, 3, 2)(1, 5, 11, 8)(6, 9, 7, 10)$.

\smallskip


\noindent{\bf Case 5.2} When ${\rm lk}(1) = C_7(10, 2, 0, 7, [x, y, z])$, it is easy to see that $(x, y, z) \in \{(5, 6, 8)$, $(5, 6, 11)$, $(5, 8, 6)$, $(5, 8, 11)$, $(5, 11, 6)$, ($5, 11, 8)$, ($8, 5, 6)$, $(8, 5, 11)$, $(8, 6, 5)$, $(8, 6, 11)$, $(8, 11, 5)$, $(8, 11, 6)$, $(9, 5, 6)$, $(9, 5, 8)$, $(9, 5, 11)$, $(9, 6, 5)$, $(9, 6, 8)$, $(9, 6, 11)$, $(9, 8, 5)$, $(9, 8, 6)$, $(9, 8, 11)$, $(11, 5, 6)$, $(11, 5, 8)$, $(11, 6, 5)$, $(11, 6, 8)$, $(11, 8, 5)$, $(11, 8, 6)\}$.

By the map $(0, 2)(3, 4)(5, 8)(6, 9)(7, 10)$ the following cases are isomorphic : $(9, 5, 8) \cong (5, 8, 6)$, $(11, 5, 8) \cong (5, 8, 11)$,
$(9, 8, 5) \cong (8, 5, 6)$, $(11, 8, 5) \cong (8, 5, 11)$, $(9, 8, 6) \cong (9, 5, 6)$, $(11, 8, 6) \cong (9, 5, 11)$ and $(11, 5, 6) \cong (9, 8, 11)$. Thus we may assume, $(x, y, z) \in \{(5, 6, 8)$, $(5, 6, 11)$, $(5, 8, 6)$, $(5, 8, 11)$, $(5, 11, 6)$, $(5, 11, 8)$, $(8, 5, 6)$, $(8, 5, 11)$, $(8, 6, 5)$, $(8, 6, 11)$, $(8, 11, 5)$, $(8, 11, 6)$, $(9, 5, 6)$, $(9, 5, 11)$, $(9, 6, 5)$, $(9, 6, 8)$, $(9, 6$, $11)$, $(9, 8, 11)$, $(11, 6, 5)$, $(11, 6, 8)\}$.

\begin{claim}\label{cl4} There does not exist a SEM for $(x, y, z) \in \{(5, 6, 8)$, $(5, 6, 11)$, $(5, 8, 6)$, $(5, 11, 6)$, $(8, 5, 6)$, $(8, 5, 11)$, $(8, 6, 5)$, $(8, 6, 11)$, $(8, 11, 6)$, $(9, 5, 6)$, $(9, 5, 11)$, $(9, 6, 5)$, $(9, 6, 8)$, $(9, 6, 11)$, $(9, 8, 11)$, $(11, 6, 5)$, $(11, 6, 8) \}$.
\end{claim}

{\sc Proof of the claim}\,\ref{cl4} When $(x, y, z) = (5, 6, 8)$ we get ${\rm lk}(8) = C_7([1, 5, 6], 9, 2, 3, 10)$. So, ${\rm lk}(10) = C_7([9, 11, 7], 3, 8, 1, 2)$ or ${\rm lk}(10) = C_7([9, 7, 11], 3, 8, 1, 2)$. In first case ${\rm lk}(7)$ has more than seven vertices. In second case ${\rm lk}(3) = C_7([4, 0, 2], 8, 10, 11, c)$ where $c \in \{5, 6, 7, 9\}$. For each values of $c$, ${\rm lk}(c)$ has more than seven vertices respectively. This is not possible. So, $(x, y, z) \neq (5, 6, 8)$. When $(x, y, z) = (5, 6, 11)$ we have ${\rm lk}(5) = C_7(7, [1, 11, 6], 0, 4, b)$ where $b\in\{8, 9, 10\}$. If $b = 8$ then ${\rm lk}(7) = C_7([8, 9, 10], 6, 0, 1, 5)$, ${\rm lk}(10) = C_7([7, 8, 9], 2, 1, 11, 6)$. This implies $9$ lies in two quadrangles. If $b = 9$ then ${\rm lk}(7) = C_7(5, 1, 0, 6, [a, b, 9])$. One can see easily that $a$ and $b$ have no suitable values in $V(K)$ which completes $K$. When $b = 10$ then we get $11$ in two quadrangles. Hence $(x, y, z) \neq (5, 6, 11)$. When $(x, y, z) = (5, 8, 6)$, ${\rm lk}(1) = C_7(10, 2, 0, 7, [5, 8, 6])$. Then considering link of $0$ we see that $56$ is both an edge and a non-edge. So, $(x, y, z) \neq (5, 8, 6)$. When $(x, y, z) = (5, 11, 6)$ then $5$ and $6$ are diagonal vertices in quadrangle in ${\rm lk}(1)$, while $56$ is an edge in ${\rm lk}(0)$. So, $(x, y, z) \neq (5, 11, 6)$. When $(x, y, z) = (8, 5, 6)$ we get ${\rm lk}(8) = C_7([1, 6, 5], 9, 2, 3, 7)$ or ${\rm lk}(8) = C_7([1, 6, 5], 3, 2, 97)$. In the first case $6$ or $3$ appear in two quadrangles. Second case implies ${\rm lk}(7) = C_7([9, 10, 11], 6, 0,1, 8)$. Then $C_5(10, 2, 8, 7,11) \subseteq {\rm lk}(9)$. Hence $(x, y, z) \neq (8, 5, 6)$.

When $(x, y, z) = (8, 5, 11)$ we have ${\rm lk}(8) = C_7([1, 11, 5], 3, 2, 9, 7)$ or ${\rm lk}(8) = C_7([1, 11, 5], 9, 2, 3$, $7)$.  If ${\rm lk}(8) = C_7([1, 11, 5], 3, 2, 9, 7)$ then ${\rm lk}(5) = C_7([8,1, 11], 4, 0, 6, 3)$ and ${\rm lk}(3) = C_7([4, 0, 2], 8, 5$, $6, a)$ for some $a \in \{7, 9, 10, 11\}$. In case $a = 7$ we get 4 in two quadrangles. When $a = 9$ or $a = 10$, then their respective links have more than seven vertices.  When $a = 11$ we have ${\rm lk}(3) = C_7([4, 0 ,2], 8, 5, 6, 11)$. Then $C_5(11, 5, 0, 2, 3)\subseteq{\rm lk}(4)$. If ${\rm lk}(8) = C_7([1, 11, 5], 9, 2, 3, 7)$ then ${\rm lk}(5) = C_7([8, 1, 11], 6, 0, 4, 9)$ or ${\rm lk}(5) = C_7([8, 1, 11], 4, 0, 6, 9)$. When ${\rm lk}(5) = C_7([8, 1, 11], 4, 0, 6, 9)$ we get ${\rm lk} (9) = C_7(10, 2, 8, 5, [6, a, b])$. We see that there are no values of $a$ and $b$ in $V(K)$ so that $K$ can be constructed. When ${\rm lk}(5) = C_7([8, 1, 11], 6, 0, 4, 9)$ we have ${\rm lk}(9) = C_7([10, 6, 7], 4, 5, 8, 2)$ or ${\rm lk}(9) = C_7([10, 7, 6], 4, 5, 8, 2)$. If ${\rm lk}(9) = C_7([10, 6, 7], 4, 5, 8, 2)$ then we get more than 7 vertices in ${\rm lk}(7)$. This is not possible. So, $(x, y, z) \neq (8, 5, 11)$.

When $(x, y, z) = (8, 6, 5)$ we get ${\rm lk}(1) = C_7([5, 6, 8], 7, 0, 2, 10)$. This implies ${\rm lk}(5) = C_7(10, [1, 8$, $6], 0, 4, a)$ where $a \in \{7, 9, 11\}$. If $a = 7$ then ${\rm lk}(7)$ has more than seven vertices. If $a = 9$ then $C_4(9, 5, 1, 2)\in{\rm lk}(10)$. If $a = 11$ we have ${\rm
lk}(10) = C_7(9, 2, 1, 5, [11, b, c])$. But $b$ and $c$ have no values in $V(K)$ such that ${\rm lk}(10)$ may be completed. Hence $(x, y, z) \neq (8, 6, 5)$. When $(x, y, z) = (8, 6, 11)$ we get ${\rm lk}(8) = C_7([6, 11, 1], 7, 3
, 2, 9)$ or ${\rm lk}(8) = C_7([6, 11, 1], 7, 9, 2, 3)$. If ${\rm lk}(8) = C_7([6, 11, 1], 7, 3, 2, 9)$ then ${\rm lk}(6) = C_7([8, 1, 11], 5, 0, 7, 9)$ or ${\rm lk}(6) = C_7([8, 1, 11], 7, 0, 5, 9)$. In the first case ${\rm lk}(7)$ has more than seven vertices. In second case $0$ or $3$ appears in two quadrangles. If ${\rm lk}(8) = C_7([6, 11, 1], 7, 9, 2, 3)$ then ${\rm lk}(6
) = C_7([8, 1, 11], 7, 0, 5, 3)$ or ${\rm lk}(6) = C_7([8, 1, 11], 5, 0, 7, 3)$. In first case ${\rm lk }(7)$ has more than seven vertices. In second case considering ${\rm lk}(7)$ we see that 9 lies in two different quadrangles. So, $(x, y, z) \neq (8, 6, 11)$.

When $(x, y, z) = (8, 11, 6)$ then ${\rm lk}(1) = C_7([6, 11, 8], 7, 0, 2, 10)$ this implies ${\rm lk}(6) = C_7([1, 8, 11]$, $7, 0, 5, 10)$ or ${\rm lk}(6) = C_7([1, 8, 11], 5, 0, 7, 10)$. In first case, considering ${\rm lk}(7)$ we see that $8$ or $11$ appear in  two quadrangles.  So, ${\rm lk}(6) = C_7([1, 8, 11], 5, 0, 7, 10)$. This implies ${\rm lk}(8) = C_7([1, 6, 11], 3$, $2, 9, 7)$ or ${\rm lk}(8) = C_7([1, 6, 11], 9, 2, 3, 7)$. In first case ${\rm lk}(7) = C_7([10, 5, 9], 8, 1, 0, 6)$. Then $9$ and $10$ are diagonal vertices of a quadrangle. This is not possible, as they are edges in ${\rm lk}(2)$. In second case, considering ${\rm lk}(7)$ we see that $3$ lies in two quadrangles. So, $(x, y, z) \neq (8, 11, 6)$. When $(x, y, z) = (9, 5, 6)$ we have ${\rm lk}(1) = C_7([6, 5, 9], 7, 0, 2, 10)$. This implies ${\rm lk}(9) = C_7([1, 6$, $5], 10, 2, 8, 7)$ or ${\rm lk}(9) = C_7([1, 6, 5], 8, 2, 10, 7)$. In the first case we get ${\rm lk}(6) = C_7(10, [1, 9, 5], 0, 7, a)$, where $a \in \{3, 4, 8, 11\}$. If $a = 3, 4$ then $a$ appears in two quadrangles. If $a = 8$ then $C_4(0, 1, 9, 8, 6)\subseteq{\rm lk}(7)$. If $a = 11$ we get $5$ in two quadrangles. Which is not allowed. If ${\rm lk}(9) = C_7([1, 6, 5], 8, 2, 10, 7)$ then ${\rm lk}(6) = C_7(10, [1, 9, 5], 0, 7, b)$ where $b \in \{3, 4, 8, 11\}$. If $b = 3, 4$, then $b$ appears in two quadrangles. If $b = 8$ then ${\rm lk}(6)$ has more than seven vertices. If $b = 11$ then 11 and 10 are diagonal vertices in quadrangle contained in ${\rm lk}(7)$, while $10\,11$ is and edge in ${\rm lk}(6)$. A contradiction. Hence $(x, y, z) \neq (9, 5, 6)$. When $(x, y, z) = (9, 5, 11)$ we have ${\rm lk}(1) = C_7([11, 5, 9], 7, 0, 2, 10)$. This implies ${\rm lk}(9) = C_7([1, 11, 5], 10, 2, 8, 7)$ or ${\rm lk}(9) = C_7([1, 11, 5], 8, 2, 10, 7)$. If ${\rm lk}(9) = C_7([1, 11, 5], 10, 2, 8, 7)$ then ${\rm lk}(5) = C_7([9, 1, 11], 6, 0, 4, 10)$ or ${\rm lk}(5) = C_7([9, 1, 11], 4, 0, 6, 10)$. In first case, considering ${\rm lk}(10)$ we get $4$ in two quadrangles. Similarly, in second case $11$ will appear in two quadrangles. If ${\rm lk}(9) = C_7([1, 11, 5], 8, 2, 10, 7)$ then ${\rm lk}(7) =
C_7([10, a, b], 6, 0, 1, 9)$ or ${\rm lk}(7) = C_7([6, a, b], 10, 9, 1, 0)$, for some $a, b \in V(K)$. In both these cases no values of $a$ and $b$ exist such that $K$ can be constructed. So, $(x, y, z) \neq (9, 5, 11)$.

When $(x, y, z) = (9, 6, 5)$ then ${\rm lk}(1) = C_7([5, 6, 9], 7, 0, 2, 10)$. So, ${\rm lk}(5) = C_7([1, 9, 6], 0 ,4, 11, 10)$, ${\rm lk}(10) = C_7([11, 8, 7], 9, 2, 1, 5)$, ${\rm lk}(7) = C_7([10, 11, 8], 6, 0, 1, 9)$, ${\rm lk}(9) = C_7([1, 5, 6], 8, 2, 10, 7)$. Then $C_6(8, 9, 1, 5, 0, 7)\subseteq{\rm lk}(6)$. Hence $(x, y, z) \neq (9, 6, 5)$. When $(x, y, z) = (9, 6, 8)$ then ${\rm lk}(1) = C_7([8, 6, 9], 7, 0, 2, 10)$. But 89 is an edge in ${\rm lk}(2)$. So this is not possible. Therefore $(x, y, z) \neq (9, 6, 8)$. When $(x, y, z) = (9 ,6, 11)$ we have ${\rm lk}(1) = C_7([11, 6, 9], 7, 0, 2, 10)$. This implies ${\rm lk}(9) = C_7([1, 11, 6], 10, 2, 8, 7)$ or ${\rm lk}(9) = C_7([1, 11, 6], 8, 2, 10, 7)$. In the first case either $11$ or $6$ appears in two quadrangles. If ${\rm lk}(9) = C_7([1, 11, 6], 8, 2, 10, 7)$ then ${\rm lk}(10) = C_7([5, 8, 7], 9, 2, 1, 11)$ or ${\rm lk}(10) = C_7([8, 5, 7], 9, 2, 1, 11)$. If ${\rm lk}(10) = C_7([5 ,8, 7], 9, 2, 1, 11)$ then ${\rm lk}(7) = C_7([10, 5, 8], 6, 0, 1, 9)$, ${\rm lk}(6) = C_7([9, 1, 11], 5, 0, 7, 8)$, ${\rm lk}(5) = C_7([10, 7, 8], 4, 0, 6, 11)$. This implies ${\rm lk}(8)$ has more than seven vertices. If ${\rm lk}(10) = C_7([8, 5, 7], 9, 2, 1, 11)$ then ${\rm lk}(7) = C_7([10, 8, 5], 6, 0, 1, 9)$. Then $C_3(5, 0, 7) \in {\rm lk}(7)$. This is not possible. So, $(x, y, z) \neq (9 ,6, 11)$. When $(x, y, z) = (9, 8, 11)$ we have ${\rm lk}(1) = C_7([11, 8, 9], 7, 0, 2, 10)$. This implies ${\rm lk}(9) = C_7([1, 11, 8]$, $2, 10, b, 7)$. It is easy to see that $b \in \{3, 4, 5, 6\}$. If $b = 3, 4$, then considering ${\rm lk}(10)$ we see that vertex $b$ or $11$ appear in two quadrangles. If $b = 5$ then ${\rm lk}(5) = C_7(0, 4, 10, 9, [7, a, 6])$ or ${\rm lk}(5) = C_7(0, 4, 7, 9, [10, a, 6])$, for some $a \in V(K)$. But we see that $a$ has no value in $V(K)$ such that $K$ can be constructed. Similarly, if $b = 6$ considering ${\rm lk}(6)$ we see that $K$ can not be constructed. Hence $(x, y, z) \neq (9, 8, 11)$.

When $(x, y ,z) = (11, 6, 5)$ we get ${\rm lk}(1) = C_7([5, 6, 11], 7, 0, 2, 10)$. Thus ${\rm lk}(5) = C_7(10, [1, 11, 6]$, $0, 4, a)$ where $a \in \{7, 8, 9\}$. If $a = 7$ then ${\rm lk}(7)$ has more than seven vertices. If $a = 8$ then ${\rm lk}(8) = C_7([10, 7, 9], 2, 3, 4, 5)$. This implies $9\,10$ is a non-edge whereas it is an edge in ${\rm lk}(2)$. If $a = 9$ then ${\rm lk}(9) = C_7(4, 5, 10, 2, [8, b, c])$ for some $b, c \in V(K)$. But $b, c$ have no values in $V(K)$ such that $K$ can be completed. So, $(x, y, z) \neq (11, 6, 5)$. When $(x, y, z) = (11, 6, 8)$ we get ${\rm lk}(1) = C_7([8, 6, 11], 7, 0, 2, 10)$, ${\rm lk}(8) = C_7([1, 11, 6], 9, 2, 3, 10)$. This implies${\rm lk}(10) = C_7([9, 5, 7], 3, 8, 1, 2)$ or ${\rm lk}(10) = C_7([9, 7, 5], 3, 8, 1, 2)$. If ${\rm lk}(10) = C_7([9, 5, 7]$, $3, 8, 1, 2)$ then ${\rm lk}(3) = C_7([4, 0, 2], 8, 10, 7, b)$ where $b \in \{5, 6, 9, 11\}$. If $b = 5$ then considering ${\rm lk}(5)$ we see that $9$ appears in two quadrangles.If $b = 6$ or $9$ then ${\rm lk}(b)$ has more than seven vertices. If $b = 11$ then ${\rm lk}(7)$ has more than seven vertices. So, ${\rm lk}(10) = C_7([9, 7, 5], 3, 8, 1, 2)$. Then ${\rm lk}(9) = C_7([7, 5, 10], 2, 8, 6, b)$ where $b \in \{4, 11\}$. If $b = 4$ then ${\rm lk}(7)$ has more than seven vertices. If $b = 11$ then ${\rm lk}(7) = C_7([5, 10, 9],11, 1, 0, 6)$. Then $C_3(5, 0, 7)\subseteq {\rm lk}(6)$. This is not possible. So, $(x, y, z) \neq (11, 6, 8)$. This completes the proof of Claim \ref{cl4}.

\smallskip

So, we have $(x, y, z) = (5, 8, 11), (5, 11, 8)$ or $(8, 11, 5)$. When $(x, y, z) = (5, 8, 11)$ we get ${\rm lk}(8) = C_7([11, 1, 5], 3, 2, 9, b)$ or ${\rm lk}(8) = C_7([5, 1, 11], 9, 2, 3, b)$ for some $b \in V(K)$.
In the first case, $i.e.$ when ${\rm lk}(8) = C_7([11, 1, 5], 3, 2, 9, b)$, we see that $b\in\{6, 7\}$. If $b = 6$ or $7$ then considering ${\rm lk}(8)$ and ${\rm lk}(6)$ we see that $5\,11$ form both - and edge and a non-edge.
This is a contradiction. If ${\rm lk}(8) = C_7([5, 1, 11], 9, 2, 3, b)$ then we get $b \in \{6, 7, 10\}$. When $b = 7$, considering ${\rm lk}(7)$ we get $3$ in two quadrangles. When $b = 10$ then ${\rm lk}(6)$ has more than seven vertices. If $b = 6$ then we get ${\rm lk}(6) = C_7([9, 10, 7], 0, 5, 8, 3)$ or ${\rm lk}(6) = C_7([10, 9, 7], 0, 5, 8, 3)$. When ${\rm lk}(6) = C_7([9, 10, 7], 0, 5, 8, 3)$ completing successively we get ${\rm lk}(10) = C_7([9, 6, 7], 4, 11, 1, 2)$, ${\rm lk}(9) = C_7([6, 7, 10], 2, 8, 11, 3)$, ${\rm lk}(3) = C_7([2, 0, 4], 11, 9, 6, 8)$, ${\rm lk}(11) = C_7([8, 5, 1], 10, 4, 3, 9)$, ${\rm lk}(7) = C_7([10, 9, 6], 0, 1, 5, 4)$, ${\rm lk}(4) = C_7([3, 2, 0], 5, 7, 10, 11)$, ${\rm lk}(5) = C_7([1, 11, 8], 6, 0, 4, 7)$. This is isomorphic to $K_2$ by the map $(0, 8, 10)(1, 6, 3, 11, 7, 2, 5, 9, 4)$. When ${\rm lk}(6) = C_7([10, 9, 7], 0, 5, 8, 3)$ completing successively we get $ {\rm lk}(10) = C_7([9, 7, 6], 3, 11, 1, 2)$, ${\rm lk}(9) = C_7([10, 6, 7], 4, 11, 8, 2)$, ${\rm lk}(4) = C_7([3, 2, 0], 5, 7, 9, 11)$, ${\rm lk}(11) = C_7([8, 5, 1], 10, 3, 4, 9)$, ${\rm lk}(3) = C_7([2, 0, 4], 11, 10, 6, 8)$, ${\rm lk}(5) = C_7([1, 11, 8], 6, 0, 4, 7)$, ${\rm lk}(7) = C_7([9, 10, 6], 0, 1, 5, 4)$. This is isomorphic to $K_3$ by the map $(0, 10, 8)(1, 3, 6, 5, 4, 9)(2, 7, 11)$.

When $(x, y, z) = (5, 11, 8)$, then ${\rm lk}(8) = C_7([1, 5, 11], 9, 2, 3, 10)$, ${\rm lk}(5) = C_7([1, 8, 11], 6, 0, 4, 7)$ this implies ${\rm lk}(7) = C_7([6, 9, 10], 4, 5, 1, 0)$ or ${\rm lk}(7) = C_7([6, 10, 9], 4, 5, 1, 0)$.
If ${\rm lk}(7) = C_7([6, 9, 10]$, $4, 5, 1, 0)$ then ${\rm lk}(10)$ has more than seven vertices. If ${\rm lk}(7) = C_7([6, 10, 9], 4, 5, 1, 0)$ then completing successively we get ${\rm lk}(10) = C_7([6, 7, 9], 2, 1, 8, 3)$, ${\rm lk}(3) = C_7([2, 0, 4], 11, 6, 10, 8)$, ${\rm lk}(6) = C_7([10, 9, 7], 0, 5, 11, 3)$, ${\rm lk}(11) = C_7([5, 1, 8], 9, 4, 3, 6)$, ${\rm lk}(4) = C_7([3, 2, 0], 5, 7, 9, 11)$, ${\rm lk}(9) = C_7([10$, $6, 7], 4, 11, 8, 2)$. It is isomorphic to $K_2$ by the map $(0, 9, 5, 2, 10, 11)(1, 3, 7, 8, 4, 6)$.

If ${\rm lk}(9) = C_7([10, 7, 6], 4, 5, 8, 2)$ then completing successively we get ${\rm lk}(6) = C_7([9, 10, 7], 0, 5$, $11, 4)$, ${\rm lk}(7) = C_7([6, 9, 10], 3, 8, 1, 0)$, ${\rm lk}(10) = C_7([7, 6, 9], 2, 1, 11, 3)$, ${\rm lk}(11) = C_7([5, 8, 1], 10, 3, 4$, $6)$, ${\rm lk}(4) = C_7([0, 2, 3], 11, 6, 9, 5)$, ${\rm lk}(3) = C_7([4, 0, 2], 8, 7, 10, 11)$. It is isomorphic to $K_3$ by the map $(0, 3)(1, 8)(5, 11)(6, 10)$.

When $(x, y, z)$ = (8, 11, 5) then ${\rm lk}(8) = C_7([1, 5, 11], 3, 2,
9, 7)$ or ${\rm lk}(8) = C_7([1, 5, 11], 9, 2, 3, 7)$. If ${\rm lk}(8) = C_7([1, 5, 11], 3, 2, 9, 7)$ then ${\rm lk}(5)  = C_7([1, 8, 11], 6, 0, 4, 10)$ or ${\rm lk}(5) = C_7([1, 8, 11], 4, 0$, $6, 10)$. If ${\rm lk}(5) = C_7([1, 8, 11], 6, 0, 4, 10)$ then ${\rm lk}(10) = C_7([6, 7, 9], 2, 1, 5, 4)$. Then ${\rm lk}(7)$ has just four faces. Which is not possible. When ${\rm lk}(5) = C_7([1, 8, 11], 4, 0, 6, 10)$ we get ${\rm lk}(10) = C_7([7, 8, 9], 2, 1, 5, 6)$. Then $C_4(2, 8, 7, 10)\subseteq{\rm lk}(9)$. This is not possible. If ${\rm lk}(8) = C_7([1, 5, 11], 9, 2$, $3, 7)$ then ${\rm lk}(5)  = C_7([1, 8, 11], 6, 0, 4, 10)$ or ${\rm lk}(5) = C_7([1, 8$, $11], 4, 0, 6, 10)$. When ${\rm lk}(5) = C_7([1, 8, 11], 4, 0, 6, 10)$ this implies ${\rm lk}(10) = C_7([6, a, b], 9, 2, 1, 5)$ or ${\rm lk}(10) = C_7([9, a, b], 6, 5, 1, 2)$, for some $a, b \in V(K)$. In both cases no values of $a$ and $b$ exists such that $K$ can be constructed. If ${\rm lk}(5) = C_7([1, 8, 11], 6, 0, 4, 10)$ then ${\rm lk}(10) = C_7([7, 6, 9], 2, 1, 5, 4)$ or ${\rm lk}(10) = C_7([9, 7, 6], 4, 5, 1, 2)$. In first case ${\rm lk}$(7) has more than seven vertices. In second case, completing successively we get ${\rm lk}(9) = C_7([10, 6, 7], 3, 11, 8, 2)$, ${\rm lk}(7) =  C_7([9, 10, 6], 0, 1, 8, 3)$, ${\rm lk}(3) = C_7([2, 0, 4], 11, 9, 7, 8)$, ${\rm lk}(11) = C_7([8, 1, 5], 6, 4, 3, 9)$ and ${\rm lk}(4) = C_7([0, 2, 3], 11, 6, 10, 5)$. It is isomorphic to $K_2$ by the map $(0, 6, 8, 4, 9, 11, 3, 10, 1)(2, 7, 1)$.


\hfill$\Box$

\section{Acknowledgement}

Part of this work was done when the first author was visiting Department of Mathematics, Indian Institute of Science during June - July 2010. We would like to thank Prof. B. Datta for numerous suggestions which led to significant improvements in the article. We would also like to thank Prof. S. C. Gupta whose suggestions proved valuable.


\smallskip

\end{document}